%
\newskip\ttglue
\font\fiverm=cmr5
\font\fivei=cmmi5
\font\fivesy=cmsy5
\font\fivebf=cmbx5
\font\sixrm=cmr6
\font\sixi=cmmi6
\font\sixsy=cmsy6
\font\sixbf=cmbx6
\font\sevenrm=cmr7
\font\eightrm=cmr8
\font\eighti=cmmi8
\font\eightsy=cmsy8
\font\eightit=cmti8
\font\eightsl=cmsl8
\font\eighttt=cmtt8
\font\eightbf=cmbx8
\font\ninerm=cmr9
\font\ninei=cmmi9
\font\ninesy=cmsy9
\font\nineit=cmti9
\font\ninesl=cmsl9
\font\ninett=cmtt9
\font\ninebf=cmbx9
\font\twelverm=cmr12
\font\twelvei=cmmi12
\font\twelvesy=cmsy12
\font\twelveit=cmti12
\font\twelvesl=cmsl12
\font\twelvett=cmtt12
\font\twelvebf=cmbx12


\def\eightpoint{\def\rm{\fam0\eightrm}  
  \textfont0=\eightrm \scriptfont0=\sixrm \scriptscriptfont0=\fiverm
  \textfont1=\eighti  \scriptfont1=\sixi  \scriptscriptfont1=\fivei
  \textfont2=\eightsy  \scriptfont2=\sixsy  \scriptscriptfont2=\fivesy
  \textfont3=\tenex  \scriptfont3=\tenex  \scriptscriptfont3=\tenex
  \textfont\itfam=\eightit  \def\it{\fam\itfam\eightit}
  \textfont\slfam=\eightsl  \def\sl{\fam\slfam\eightsl}
  \textfont\ttfam=\eighttt  \def\tt{\fam\ttfam\eighttt}
  \textfont\bffam=\eightbf  \scriptfont\bffam=\sixbf
    \scriptscriptfont\bffam=\fivebf  \def\bf{\fam\bffam\eightbf}
  \tt  \ttglue=.5em plus.25em minus.15em
  \normalbaselineskip=9pt
  \setbox\strutbox=\hbox{\vrule height7pt depth2pt width0pt}
  \let\sc=\sixrm  \let\big=\eightbig \normalbaselines\rm}

\def\eightbig#1{{\hbox{$\textfont0=\ninerm\textfont2=\ninesy
        \left#1\vbox to6.5pt{}\right.$}}}


\def\ninepoint{\def\rm{\fam0\ninerm}  
  \textfont0=\ninerm \scriptfont0=\sixrm \scriptscriptfont0=\fiverm
  \textfont1=\ninei  \scriptfont1=\sixi  \scriptscriptfont1=\fivei
  \textfont2=\ninesy  \scriptfont2=\sixsy  \scriptscriptfont2=\fivesy
  \textfont3=\tenex  \scriptfont3=\tenex  \scriptscriptfont3=\tenex
  \textfont\itfam=\nineit  \def\it{\fam\itfam\nineit}
  \textfont\slfam=\ninesl  \def\sl{\fam\slfam\ninesl}
  \textfont\ttfam=\ninett  \def\tt{\fam\ttfam\ninett}
  \textfont\bffam=\ninebf  \scriptfont\bffam=\sixbf
    \scriptscriptfont\bffam=\fivebf  \def\bf{\fam\bffam\ninebf}
  \tt  \ttglue=.5em plus.25em minus.15em
  \normalbaselineskip=11pt
  \setbox\strutbox=\hbox{\vrule height8pt depth3pt width0pt}
  \let\sc=\sevenrm  \let\big=\ninebig \normalbaselines\rm}

\def\ninebig#1{{\hbox{$\textfont0=\tenrm\textfont2=\tensy
        \left#1\vbox to7.25pt{}\right.$}}}


\def\twelvepoint{\def\rm{\fam0\twelverm}  
  \textfont0=\twelverm \scriptfont0=\eightrm \scriptscriptfont0=\sixrm
  \textfont1=\twelvei  \scriptfont1=\eighti  \scriptscriptfont1=\sixi
  \textfont2=\twelvesy  \scriptfont2=\eightsy  \scriptscriptfont2=\sixsy
  \textfont3=\tenex  \scriptfont3=\tenex  \scriptscriptfont3=\tenex
  \textfont\itfam=\twelveit  \def\it{\fam\itfam\twelveit}
  \textfont\slfam=\twelvesl  \def\sl{\fam\slfam\twelvesl}
  \textfont\ttfam=\twelvett  \def\tt{\fam\ttfam\twelvett}
  \textfont\bffam=\twelvebf  \scriptfont\bffam=\eightbf
    \scriptscriptfont\bffam=\sixbf  \def\bf{\fam\bffam\twelvebf}
  \tt  \ttglue=.5em plus.25em minus.15em
  \normalbaselineskip=11pt
  \setbox\strutbox=\hbox{\vrule height8pt depth3pt width0pt}
  \let\sc=\sevenrm  \let\big=\twelvebig \normalbaselines\rm}

\def\twelvebig#1{{\hbox{$\textfont0=\tenrm\textfont2=\tensy
        \left#1\vbox to7.25pt{}\right.$}}}
\magnification=\magstep1
\def\firstpage{1}
\pageno=\firstpage
\font\fiverm=cmr5
\font\sevenrm=cmr7
\font\sevenbf=cmbx7
\font\eightrm=cmr8
\font\eightbf=cmbx8
\font\ninerm=cmr9
\font\ninebf=cmbx9
\font\tenbf=cmbx10
\font\magtenbf=cmbx10 scaled\magstep1

\font\nineeufm=eufm9
\font\magnineeufm=eufm9 scaled\magstep1

%
%
\newskip\ttglue
\font\fiverm=cmr5
\font\fivei=cmmi5
\font\fivesy=cmsy5
\font\fivebf=cmbx5
\font\sixrm=cmr6
\font\sixi=cmmi6
\font\sixsy=cmsy6
\font\sixbf=cmbx6
\font\sevenrm=cmr7
\font\eightrm=cmr8
\font\eighti=cmmi8
\font\eightsy=cmsy8
\font\eightit=cmti8
\font\eightsl=cmsl8
\font\eighttt=cmtt8
\font\eightbf=cmbx8
\font\ninerm=cmr9
\font\ninei=cmmi9
\font\ninesy=cmsy9
\font\nineit=cmti9
\font\ninesl=cmsl9
\font\ninett=cmtt9
\font\ninebf=cmbx9
\font\twelverm=cmr12
\font\twelvei=cmmi12
\font\twelvesy=cmsy12
\font\twelveit=cmti12
\font\twelvesl=cmsl12
\font\twelvett=cmtt12
\font\twelvebf=cmbx12


\def\eightpoint{\def\rm{\fam0\eightrm}  
  \textfont0=\eightrm \scriptfont0=\sixrm \scriptscriptfont0=\fiverm
  \textfont1=\eighti  \scriptfont1=\sixi  \scriptscriptfont1=\fivei
  \textfont2=\eightsy  \scriptfont2=\sixsy  \scriptscriptfont2=\fivesy
  \textfont3=\tenex  \scriptfont3=\tenex  \scriptscriptfont3=\tenex
  \textfont\itfam=\eightit  \def\it{\fam\itfam\eightit}
  \textfont\slfam=\eightsl  \def\sl{\fam\slfam\eightsl}
  \textfont\ttfam=\eighttt  \def\tt{\fam\ttfam\eighttt}
  \textfont\bffam=\eightbf  \scriptfont\bffam=\sixbf
    \scriptscriptfont\bffam=\fivebf  \def\bf{\fam\bffam\eightbf}
  \tt  \ttglue=.5em plus.25em minus.15em
  \normalbaselineskip=9pt
  \setbox\strutbox=\hbox{\vrule height7pt depth2pt width0pt}
  \let\sc=\sixrm  \let\big=\eightbig \normalbaselines\rm}

\def\eightbig#1{{\hbox{$\textfont0=\ninerm\textfont2=\ninesy
        \left#1\vbox to6.5pt{}\right.$}}}


\def\ninepoint{\def\rm{\fam0\ninerm}  
  \textfont0=\ninerm \scriptfont0=\sixrm \scriptscriptfont0=\fiverm
  \textfont1=\ninei  \scriptfont1=\sixi  \scriptscriptfont1=\fivei
  \textfont2=\ninesy  \scriptfont2=\sixsy  \scriptscriptfont2=\fivesy
  \textfont3=\tenex  \scriptfont3=\tenex  \scriptscriptfont3=\tenex
  \textfont\itfam=\nineit  \def\it{\fam\itfam\nineit}
  \textfont\slfam=\ninesl  \def\sl{\fam\slfam\ninesl}
  \textfont\ttfam=\ninett  \def\tt{\fam\ttfam\ninett}
  \textfont\bffam=\ninebf  \scriptfont\bffam=\sixbf
    \scriptscriptfont\bffam=\fivebf  \def\bf{\fam\bffam\ninebf}
  \tt  \ttglue=.5em plus.25em minus.15em
  \normalbaselineskip=11pt
  \setbox\strutbox=\hbox{\vrule height8pt depth3pt width0pt}
  \let\sc=\sevenrm  \let\big=\ninebig \normalbaselines\rm}

\def\ninebig#1{{\hbox{$\textfont0=\tenrm\textfont2=\tensy
        \left#1\vbox to7.25pt{}\right.$}}}


\def\twelvepoint{\def\rm{\fam0\twelverm}  
  \textfont0=\twelverm \scriptfont0=\eightrm \scriptscriptfont0=\sixrm
  \textfont1=\twelvei  \scriptfont1=\eighti  \scriptscriptfont1=\sixi
  \textfont2=\twelvesy  \scriptfont2=\eightsy  \scriptscriptfont2=\sixsy
  \textfont3=\tenex  \scriptfont3=\tenex  \scriptscriptfont3=\tenex
  \textfont\itfam=\twelveit  \def\it{\fam\itfam\twelveit}
  \textfont\slfam=\twelvesl  \def\sl{\fam\slfam\twelvesl}
  \textfont\ttfam=\twelvett  \def\tt{\fam\ttfam\twelvett}
  \textfont\bffam=\twelvebf  \scriptfont\bffam=\eightbf
    \scriptscriptfont\bffam=\sixbf  \def\bf{\fam\bffam\twelvebf}
  \tt  \ttglue=.5em plus.25em minus.15em
  \normalbaselineskip=11pt
  \setbox\strutbox=\hbox{\vrule height8pt depth3pt width0pt}
  \let\sc=\sevenrm  \let\big=\twelvebig \normalbaselines\rm}

\def\twelvebig#1{{\hbox{$\textfont0=\tenrm\textfont2=\tensy
        \left#1\vbox to7.25pt{}\right.$}}}
\catcode`\@=11
%

\def\undefine#1{\let#1\undefined}
\def\newsymbol#1#2#3#4#5{\let\next@\relax
 \ifnum#2=\@ne\let\next@\msafam@\else
 \ifnum#2=\tw@\let\next@\msbfam@\fi\fi
 \mathchardef#1="#3\next@#4#5}
\def\mathhexbox@#1#2#3{\relax
 \ifmmode\mathpalette{}{\m@th\mathchar"#1#2#3}%
 \else\leavevmode\hbox{$\m@th\mathchar"#1#2#3$}\fi}
\def\hexnumber@#1{\ifcase#1 0\or 1\or 2\or 3\or 4\or 5\or 6\or 7\or 8\or
 9\or A\or B\or C\or D\or E\or F\fi}

\font\tenmsa=msam10
\font\sevenmsa=msam7
\font\fivemsa=msam5
\newfam\msafam
\textfont\msafam=\tenmsa
\scriptfont\msafam=\sevenmsa
\scriptscriptfont\msafam=\fivemsa
\edef\msafam@{\hexnumber@\msafam}
\mathchardef\dabar@"0\msafam@39
\def\dashrightarrow{\mathrel{\dabar@\dabar@\mathchar"0\msafam@4B}}
\def\dashleftarrow{\mathrel{\mathchar"0\msafam@4C\dabar@\dabar@}}

\def\ulcorner{\delimiter"4\msafam@70\msafam@70 }
\def\urcorner{\delimiter"5\msafam@71\msafam@71 }
\def\llcorner{\delimiter"4\msafam@78\msafam@78 }
\def\lrcorner{\delimiter"5\msafam@79\msafam@79 }
\def\yen{{\mathhexbox@\msafam@55}}
\def\checkmark{{\mathhexbox@\msafam@58}}
\def\circledR{{\mathhexbox@\msafam@72}}
\def\maltese{{\mathhexbox@\msafam@7A}}

\font\tenmsb=msbm10
\font\sevenmsb=msbm7
\font\fivemsb=msbm5
\newfam\msbfam
\textfont\msbfam=\tenmsb
\scriptfont\msbfam=\sevenmsb
\scriptscriptfont\msbfam=\fivemsb
\edef\msbfam@{\hexnumber@\msbfam}
\def\Bbb#1{{\fam\msbfam\relax#1}}
\def\widehat#1{\setbox\z@\hbox{$\m@th#1$}%
 \ifdim\wd\z@>\tw@ em\mathaccent"0\msbfam@5B{#1}%
 \else\mathaccent"0362{#1}\fi}
\def\widetilde#1{\setbox\z@\hbox{$\m@th#1$}%
 \ifdim\wd\z@>\tw@ em\mathaccent"0\msbfam@5D{#1}%
 \else\mathaccent"0365{#1}\fi}
\font\teneufm=eufm10
\font\seveneufm=eufm7
\font\fiveeufm=eufm5
\newfam\eufmfam
\textfont\eufmfam=\teneufm
\scriptfont\eufmfam=\seveneufm
\scriptscriptfont\eufmfam=\fiveeufm

\catcode`\@=11
\newsymbol\boxdot 1200
\newsymbol\boxplus 1201
\newsymbol\boxtimes 1202
\newsymbol\square 1003
\newsymbol\blacksquare 1004
\newsymbol\centerdot 1205
\newsymbol\lozenge 1006
\newsymbol\blacklozenge 1007
\newsymbol\circlearrowright 1308
\newsymbol\circlearrowleft 1309
\undefine\rightleftharpoons
\newsymbol\rightleftharpoons 130A
\newsymbol\leftrightharpoons 130B
\newsymbol\boxminus 120C
\newsymbol\Vdash 130D
\newsymbol\Vvdash 130E
\newsymbol\vDash 130F
\newsymbol\twoheadrightarrow 1310
\newsymbol\twoheadleftarrow 1311
\newsymbol\leftleftarrows 1312
\newsymbol\rightrightarrows 1313
\newsymbol\upuparrows 1314
\newsymbol\downdownarrows 1315
\newsymbol\upharpoonright 1316
 
\newsymbol\downharpoonright 1317
\newsymbol\upharpoonleft 1318
\newsymbol\downharpoonleft 1319
\newsymbol\rightarrowtail 131A
\newsymbol\leftarrowtail 131B
\newsymbol\leftrightarrows 131C
\newsymbol\rightleftarrows 131D
\newsymbol\Lsh 131E
\newsymbol\Rsh 131F
\newsymbol\rightsquigarrow 1320
\newsymbol\leftrightsquigarrow 1321
\newsymbol\looparrowleft 1322
\newsymbol\looparrowright 1323
\newsymbol\circeq 1324
\newsymbol\succsim 1325
\newsymbol\gtrsim 1326
\newsymbol\gtrapprox 1327
\newsymbol\multimap 1328
\newsymbol\therefore 1329
\newsymbol\because 132A
\newsymbol\doteqdot 132B
 
\newsymbol\triangleq 132C
\newsymbol\precsim 132D
\newsymbol\lesssim 132E
\newsymbol\lessapprox 132F
\newsymbol\eqslantless 1330
\newsymbol\eqslantgtr 1331
\newsymbol\curlyeqprec 1332
\newsymbol\curlyeqsucc 1333
\newsymbol\preccurlyeq 1334
\newsymbol\leqq 1335
\newsymbol\leqslant 1336
\newsymbol\lessgtr 1337
\newsymbol\backprime 1038
\newsymbol\risingdotseq 133A
\newsymbol\fallingdotseq 133B
\newsymbol\succcurlyeq 133C
\newsymbol\geqq 133D
\newsymbol\geqslant 133E
\newsymbol\gtrless 133F
\newsymbol\sqsubset 1340
\newsymbol\sqsupset 1341
\newsymbol\vartriangleright 1342
\newsymbol\vartriangleleft 1343
\newsymbol\trianglerighteq 1344
\newsymbol\trianglelefteq 1345
\newsymbol\bigstar 1046
\newsymbol\between 1347
\newsymbol\blacktriangledown 1048
\newsymbol\blacktriangleright 1349
\newsymbol\blacktriangleleft 134A
\newsymbol\vartriangle 134D
\newsymbol\blacktriangle 104E
\newsymbol\triangledown 104F
\newsymbol\eqcirc 1350
\newsymbol\lesseqgtr 1351
\newsymbol\gtreqless 1352
\newsymbol\lesseqqgtr 1353
\newsymbol\gtreqqless 1354
\newsymbol\Rrightarrow 1356
\newsymbol\Lleftarrow 1357
\newsymbol\veebar 1259
\newsymbol\barwedge 125A
\newsymbol\doublebarwedge 125B
\undefine\angle
\newsymbol\angle 105C
\newsymbol\measuredangle 105D
\newsymbol\sphericalangle 105E
\newsymbol\varpropto 135F
\newsymbol\smallsmile 1360
\newsymbol\smallfrown 1361
\newsymbol\Subset 1362
\newsymbol\Supset 1363
\newsymbol\Cup 1264
 
\newsymbol\Cap 1265
 
\newsymbol\curlywedge 1266
\newsymbol\curlyvee 1267
\newsymbol\leftthreetimes 1268
\newsymbol\rightthreetimes 1269
\newsymbol\subseteqq 136A
\newsymbol\supseteqq 136B
\newsymbol\bumpeq 136C
\newsymbol\Bumpeq 136D
\newsymbol\lll 136E
 
\newsymbol\ggg 136F
 
\newsymbol\circledS 1073
\newsymbol\pitchfork 1374
\newsymbol\dotplus 1275
\newsymbol\backsim 1376
\newsymbol\backsimeq 1377
\newsymbol\complement 107B
\newsymbol\intercal 127C
\newsymbol\circledcirc 127D
\newsymbol\circledast 127E
\newsymbol\circleddash 127F
\newsymbol\lvertneqq 2300
\newsymbol\gvertneqq 2301
\newsymbol\nleq 2302
\newsymbol\ngeq 2303
\newsymbol\nless 2304
\newsymbol\ngtr 2305
\newsymbol\nprec 2306
\newsymbol\nsucc 2307
\newsymbol\lneqq 2308
\newsymbol\gneqq 2309
\newsymbol\nleqslant 230A
\newsymbol\ngeqslant 230B
\newsymbol\lneq 230C
\newsymbol\gneq 230D
\newsymbol\npreceq 230E
\newsymbol\nsucceq 230F
\newsymbol\precnsim 2310
\newsymbol\succnsim 2311
\newsymbol\lnsim 2312
\newsymbol\gnsim 2313
\newsymbol\nleqq 2314
\newsymbol\ngeqq 2315
\newsymbol\precneqq 2316
\newsymbol\succneqq 2317
\newsymbol\precnapprox 2318
\newsymbol\succnapprox 2319
\newsymbol\lnapprox 231A
\newsymbol\gnapprox 231B
\newsymbol\nsim 231C
\newsymbol\ncong 231D
\newsymbol\diagup 201E
\newsymbol\diagdown 201F
\newsymbol\varsubsetneq 2320
\newsymbol\varsupsetneq 2321
\newsymbol\nsubseteqq 2322
\newsymbol\nsupseteqq 2323
\newsymbol\subsetneqq 2324
\newsymbol\supsetneqq 2325
\newsymbol\varsubsetneqq 2326
\newsymbol\varsupsetneqq 2327
\newsymbol\subsetneq 2328
\newsymbol\supsetneq 2329
\newsymbol\nsubseteq 232A
\newsymbol\nsupseteq 232B
\newsymbol\nparallel 232C
\newsymbol\nmid 232D
\newsymbol\nshortmid 232E
\newsymbol\nshortparallel 232F
\newsymbol\nvdash 2330
\newsymbol\nVdash 2331
\newsymbol\nvDash 2332
\newsymbol\nVDash 2333
\newsymbol\ntrianglerighteq 2334
\newsymbol\ntrianglelefteq 2335
\newsymbol\ntriangleleft 2336
\newsymbol\ntriangleright 2337
\newsymbol\nleftarrow 2338
\newsymbol\nrightarrow 2339
\newsymbol\nLeftarrow 233A
\newsymbol\nRightarrow 233B
\newsymbol\nLeftrightarrow 233C
\newsymbol\nleftrightarrow 233D
\newsymbol\divideontimes 223E
\newsymbol\varnothing 203F
\newsymbol\nexists 2040
\newsymbol\Finv 2060
\newsymbol\Game 2061
\newsymbol\mho 2066
\newsymbol\eth 2067
\newsymbol\eqsim 2368
\newsymbol\beth 2069
\newsymbol\gimel 206A
\newsymbol\daleth 206B
\newsymbol\lessdot 236C
\newsymbol\gtrdot 236D
\newsymbol\ltimes 226E
\newsymbol\rtimes 226F
\newsymbol\shortmid 2370
\newsymbol\shortparallel 2371
\newsymbol\smallsetminus 2272
\newsymbol\thicksim 2373
\newsymbol\thickapprox 2374
\newsymbol\approxeq 2375
\newsymbol\succapprox 2376
\newsymbol\precapprox 2377
\newsymbol\curvearrowleft 2378
\newsymbol\curvearrowright 2379
\newsymbol\digamma 207A
\newsymbol\varkappa 207B
\newsymbol\Bbbk 207C
\newsymbol\hslash 207D
\undefine\hbar
\newsymbol\hbar 207E
\newsymbol\backepsilon 237F

%
\newcount\marknumber	\marknumber=1
\newcount\countdp \newcount\countwd \newcount\countht 
%
%
\ifx\pdfoutput\undefined
\def\rgboo#1{}
\def\postscript#1{\special{" #1}}		
\postscript{
	/bd {bind def} bind def
	/fsd {findfont exch scalefont def} bd
	/sms {setfont moveto show} bd
	/ms {moveto show} bd
	/pdfmark where		
	{pop} {userdict /pdfmark /cleartomark load put} ifelse
	[ /PageMode /UseOutlines		
	/DOCVIEW pdfmark}
\def\bookmark#1#2{\postscript{		
	[ /Dest /MyDest\the\marknumber /View [ /XYZ null null null ] /DEST pdfmark
	[ /Title (#2) /Count #1 /Dest /MyDest\the\marknumber /OUT pdfmark}%
	\advance\marknumber by1}
\def\pdfclink#1#2#3{%
	\hskip-.25em\setbox0=\hbox{#2}%
		\countdp=\dp0 \countwd=\wd0 \countht=\ht0%
		\divide\countdp by65536 \divide\countwd by65536%
			\divide\countht by65536%
		\advance\countdp by1 \advance\countwd by1%
			\advance\countht by1%
		\def\linkdp{\the\countdp} \def\linkwd{\the\countwd}%
			\def\linkht{\the\countht}%
	\postscript{
		[ /Rect [ -1.5 -\linkdp.0 0\linkwd.0 0\linkht.5 ] 
		/Border [ 0 0 0 ]
		/Action << /Subtype /URI /URI (#3) >>
		/Subtype /Link
		/ANN pdfmark}{\rgb{#1}{#2}}}
%
%
\else
\def\rgboo#1{\pdfliteral{#1 rg #1 RG}}
\pdfcatalog{/PageMode /UseOutlines}		
\def\bookmark#1#2{
	\pdfdest num \marknumber xyz
	\pdfoutline goto num \marknumber count #1 {#2}
	\advance\marknumber by1}
\def\pdfklink#1#2{%
	\noindent\pdfstartlink user
		{/Subtype /Link
		/Border [ 0 0 0 ]
		/A << /S /URI /URI (#2) >>}{\rgb{1 0 0}{#1}}%
	\pdfendlink}
\fi

\def\rgbo#1#2{\rgboo{#1}#2\rgboo{0 0 0}}
\def\rgb#1#2{\mark{#1}\rgbo{#1}{#2}\mark{0 0 0}}
\def\pdfklink#1#2{\pdfclink{1 0 0}{#1}{#2}}
\def\pdflink#1{\pdfklink{#1}{#1}}
%
%
\newcount\seccount  
\newcount\subcount  
\newcount\clmcount  
\newcount\equcount  
\newcount\refcount  
\newcount\demcount  
\newcount\execount  
\newcount\procount  
\seccount=0
\equcount=1
\clmcount=1
\subcount=1
\refcount=1
\demcount=0
\execount=0
\procount=0
%
\def\proof{\medskip\noindent{\bf Proof.\ }}
\def\proofof(#1){\medskip\noindent{\bf Proof of \csname c#1\endcsname.\ }}
\def\qed{\hfill{\sevenbf QED}\par\medskip}
\def\references{\bigskip\noindent\hbox{\bf References}\medskip
                \ifx\pdflink\undefined\else\bookmark{0}{References}\fi}
\def\addref#1{\expandafter\xdef\csname r#1\endcsname{\number\refcount}
    \global\advance\refcount by 1}

\def\nextremark #1\par{\item{$\circ$} #1}
\def\firstremark #1\par{\bigskip\noindent{\bf Remarks.}
     \smallskip\nextremark #1\par}
\def\abstract#1\par{{\baselineskip=10pt
    \eightpoint\narrower\noindent{\eightbf Abstract.} #1\par}}
%
\def\equtag#1{\expandafter\xdef\csname e#1\endcsname{(\number\seccount.\number\equcount)}
              \global\advance\equcount by 1}
\def\equation(#1){\equtag{#1}\eqno\csname e#1\endcsname}
\def\equ(#1){\hskip-0.03em\csname e#1\endcsname}
%
\def\clmtag#1#2{\expandafter\xdef\csname cn#2\endcsname{\number\seccount.\number\clmcount}
                \expandafter\xdef\csname c#2\endcsname{#1~\number\seccount.\number\clmcount}
                \global\advance\clmcount by 1}
\def\claim #1(#2) #3\par{\clmtag{#1}{#2}
    \vskip.1in\medbreak\noindent
    {\bf \csname c#2\endcsname .\ }{\sl #3}\par
    \ifdim\lastskip<\medskipamount
    \removelastskip\penalty55\medskip\fi}
\def\clm(#1){\csname c#1\endcsname}
\def\clmno(#1){\csname cn#1\endcsname}
%
\def\sectag#1{\global\advance\seccount by 1
              \expandafter\xdef\csname sectionname\endcsname{\number\seccount. #1}
              \equcount=1 \clmcount=1 \subcount=1 \execount=0 \procount=0}
\def\section#1\par{\vskip0pt plus.1\vsize\penalty-40
    \vskip0pt plus -.1\vsize\bigskip\bigskip
    \sectag{#1}
    \message{\sectionname}\leftline{\magtenbf\sectionname}
    \nobreak\smallskip\noindent
    \ifx\pdflink\undefined
    \else
      \bookmark{0}{\sectionname}
    \fi}
%
\def\subtag#1{\expandafter\xdef\csname subsectionname\endcsname{\number\seccount.\number\subcount. #1}
              \global\advance\subcount by 1}
\def\subsection#1\par{\vskip0pt plus.05\vsize\penalty-20
    \vskip0pt plus -.05\vsize\medskip\medskip
    \subtag{#1}
    \message{\subsectionname}\leftline{\tenbf\subsectionname}
    \nobreak\smallskip\noindent
    \ifx\pdflink\undefined
    \else
      \bookmark{0}{.... \subsectionname}  
    \fi}
%
\def\demtag#1#2{\global\advance\demcount by 1
              \expandafter\xdef\csname de#2\endcsname{#1~\number\demcount}}
\def\demo #1(#2) #3\par{
  \demtag{#1}{#2}
  \vskip.1in\medbreak\noindent
  {\bf #1 \number\demcount.\enspace}
  {\rm #3}\par
  \ifdim\lastskip<\medskipamount
  \removelastskip\penalty55\medskip\fi}
\def\dem(#1){\csname de#1\endcsname}
%
\def\exetag#1{\global\advance\execount by 1
              \expandafter\xdef\csname ex#1\endcsname{Exercise~\number\seccount.\number\execount}}
\def\exercise(#1) #2\par{
  \exetag{#1}
  \vskip.1in\medbreak\noindent
  {\bf Exercise \number\execount.}
  {\rm #2}\par
  \ifdim\lastskip<\medskipamount
  \removelastskip\penalty55\medskip\fi}
\def\exe(#1){\csname ex#1\endcsname}
%
\def\protag#1{\global\advance\procount by 1
              \expandafter\xdef\csname pr#1\endcsname{\number\seccount.\number\procount}}
\def\problem(#1) #2\par{
  \ifnum\procount=0
    \parskip=6pt
    \vbox{\bigskip\centerline{\bf Problems \number\seccount}\nobreak\medskip}
  \fi
  \protag{#1}
  \item{\number\procount.} #2}
\def\pro(#1){Problem \csname pr#1\endcsname}
%
%
%
\def\rightheadline{\hfil}
\def\leftheadline{\sevenrm\hfil HANS KOCH\hfil}
\headline={\ifnum\pageno=\firstpage\hfil\else
\ifodd\pageno{{\fiverm\rightheadline}\number\pageno}
\else{\number\pageno\fiverm\leftheadline}\fi\fi}
\footline={\ifnum\pageno=\firstpage\hss\tenrm\folio\hss\else\hss\fi}

\let\cl=\centerline

\let\eps=\varepsilon
\let\sss=\scriptscriptstyle

\def\AA{{\cal A}}
\def\BB{{\cal B}}

\def\EE{{\cal E}}
\def\FF{{\cal F}}
\def\GG{{\cal G}}
\def\HH{{\cal H}}
\def\II{{\cal I}}
\def\JJ{{\cal J}}

\def\LL{{\cal L}}

\def\RR{{\cal R}}

\def\TT{{\cal T}}

\def\rmC{\mathop{\rm C}\nolimits}
\def\rmL{\mathop{\rm L}\nolimits}
\def\id{\mathop{\rm I}\nolimits}

\def\diag{\mathop{\rm diag}\nolimits}

\def\Re{\mathop{\rm Re}\nolimits}
\def\Im{\mathop{\rm Im}\nolimits}
%
\newfam\dsfam
\def\mathds #1{{\fam\dsfam\tends #1}}

\font\tends=dsrom10
\font\eightds=dsrom8
\textfont\dsfam=\tends
\scriptfont\dsfam=\eightds
%

\def\integer{{\mathds Z}}

\def\real{{\mathds R}}
\def\complex{{\mathds C}}

\def\proj{{\Bbb P}}

\def\bskip{\bigskip\noindent}

\def\bdot{\hbox{\bf .}}
\def\bcomma{\hbox{\bf ,}}
\def\defeq{\mathrel{\mathop=^{\sss\rm def}}}
\def\half{{1\over 2}}
\def\third{{1\over 3}}
\def\quarter{{1\over 4}}
\def\thalf{{\textstyle\half}}

\def\twomat#1#2#3#4{\left[\matrix{#1&#2\cr#3&#4\cr}\right]}

%

%

%

%

\input miniltx

\ifx\pdfoutput\undefined
  \def\Gin@driver{dvips.def}  
\else
  \def\Gin@driver{pdftex.def} 
\fi
 
\input graphicx.sty
\resetatcatcode
%
\font\tenib=cmmib10

\def\bmy{{\hbox{\tenib y}}}
\let\Kt\bmk

\def\buH{{\hbox{\magnineeufm H}}}
\def\buc{{\hbox{\magnineeufm c}}}
\def\buh{{\hbox{\magnineeufm h}}}

\def\buu{{\hbox{\magnineeufm u}}}
\def\buv{{\hbox{\magnineeufm v}}}
\def\buw{{\hbox{\magnineeufm w}}}

\def\sbuv{{\hbox{\nineeufm v}}}
\font\sevenib=cmmib7
\def\sbmy{{\hbox{\sevenib y}}}
\let\ub\underbar

\let\Flow\Phi
\let\FLOW\Theta
\def\id{{\rm I}}

\def\rmS{{\rm S}}

\def\dist{{\rm dist}}
\def\range{{\rm range}}

\def\even{0}
\def\odd{1}
\def\pars{\tau}
\def\tinyskip{\hskip.7pt}
\def\ssigma{{\tinyskip\sigma}}
\def\sBB{{\scriptstyle\BB}}
\def\BBB{\BB_{r,\rho}^{\ssigma,\pars}}
\def\snone{{\sixrm none}}
\def\bdot{\hbox{\bf .}}
\def\bcomma{\hbox{\bf ,}}
\def\hdots{\line{\leaders\hbox to 0.5em{\hss .\hss}\hfil}}
\def\bskip{\bigskip\noindent}

\def\Langle{\bigl\langle}
\def\Rangle{\bigr\rangle}
\def\ltstrong{\mathrel{<\hskip-8.2pt\raise1.29pt\hbox{$\sss<$}}}
\def\gtstrong{\mathrel{\raise1.29pt\hbox{$\sss>$}\hskip-8.2pt>}}
\def\nobla{\raise.35em\rlap{\kern.32em.}\nabla}
\def\nablai{\nabla_{\!1}}
\def\nablao{\nabla_{\!0}}
\def\nablas{\nabla_{\!\sigma}}
\def\Hamiltonian{{\Bbb H}}
\def\Lagrangian{{\Bbb L}}
\def\Potential{V}
\def\circle{{\Bbb S}}
\def\stwovec#1#2{{\eightpoint\left[\matrix{#1\cr#2\cr}\right]}}
\def\stwomat#1#2#3#4{{\eightpoint\left[\matrix{#1&#2\cr#3&#4\cr}\right]}}
\def\lab#1 {\pdfclink{0 0 1}{#1}{http://web.ma.utexas.edu/users/koch/papers/breathers/anim/#1.gif}}
\addref{Kato}
\addref{AGb}
\addref{AGa}
\addref{AGT}
\addref{Rabi}
\addref{BK}
\addref{AKT}
\addref{BerIzr}
\addref{APankov}
\addref{Galla}
\addref{GFii}
\addref{FLMii}
\addref{PS}
\addref{AKii}
\addref{Files}
\addref{Ada}
\addref{Gnat}
\addref{IEEE}
\addref{MPFR}
\def\leftheadline{\sixrm\hfil ARIOLI \& KOCH\hfil}
\def\rightheadline{\sevenrm\hfil breathers and multi-breathers\hfil}
%
\cl{{\magtenbf Some breathers and multi-breathers for FPU-type chains}}
\bigskip

\cl{
Gianni Arioli
\footnote{$^1$}
{\eightpoint\hskip-2.9em
Department of Mathematics, Politecnico di Milano,
Piazza Leonardo da Vinci 32, 20133 Milano.
}
$^{\!\!\!,\!\!}$
\footnote{$^2$}
{\eightpoint\hskip-2.6em
Supported in part by the PRIN project
``Equazioni alle derivate parziali e
disuguaglianze analitico-geometriche associate''.}
and Hans Koch
\footnote{$^3$}
{\eightpoint\hskip-2.7em
Department of Mathematics, The University of Texas at Austin,
Austin, TX 78712.}
}

\bigskip
\abstract
We consider several breather solutions for FPU-type chains
that have been found numerically.
Using computer-assisted techniques,
we prove that there exist true solutions nearby,
and in some cases, we determine whether or not
the solution is spectrally stable.
Symmetry properties are considered as well.
In addition, we construct solutions that are close to
(possibly infinite) sums of breather solutions.

\section Introduction

We consider a system of interacting particles
described by the equation
$$
\omega^2\ddot q_j=\phi'(q_{j+1}-q_j)-\phi'(q_j-q_{j-1})-\psi'(q_j)\,,
\qquad j\in\integer\,,
\equation(Mainj)
$$
where $\phi(x)=\half\phi_2 x^2+\third\phi_3 x^3+\quarter\phi_4 x^4$
and $\psi(x)=\half\psi_2 x^2+\quarter\psi_4 x^4$,
with $\phi_3$ and $\phi_4$ not both zero.  
If $\phi$ and $\psi$ are given,
then the parameter $\omega$ simply fixes a time scale.

The equation \equ(Mainj) with $\psi=0$ is know as the Fermi-Pasta-Ulam (FPU) model:
the $\alpha$-model if $\phi_4=0$, or the $\beta$-model if $\phi_3=0$.
Models of this type have been studied extensively in connection
with the problem of equipartition of energy in systems
with a large number of interacting particles.
Recent surveys can be found in [\rBK,\rBerIzr,\rGalla,\rGFii].

Our goal is to construct solutions that are periodic in time,
and in some cases, to determine whether they are (spectrally) stable or not.
By choosing the value $\omega$ appropriately,
it suffices to consider solutions that are periodic with fundamental period $2\pi$.
We are interested in solutions that decrease rapidly in $|j|$,
also referred to as ``breathers'',
and in solutions that are close to sums of such breathers.

Most of the existing work on breather solutions involves
numerical computations or other types of approximations.
For simplicity, the computations often focus on solutions that have a reflection symmetry
$$
\breve q=q\,,\qquad\breve q_j(t)=\varpi q_{\sigma-j}(t-\delta t)\,,\qquad
j\in\integer\,,\quad t\in\real\,,
\equation(symmetry)
$$
where $\sigma$ is an integer, $\varpi$ is one of $\pm$, and $\delta t\in\{0,\pi\}$.
Breathers that have such a symmetry
are commonly referred to as being site-centered if $\sigma$ is even,
or bond-centered if $\sigma$ is odd.
Due to the translation-invariance of the equation \equ(Mainj),
it suffices to consider $\sigma\in\{0,1\}$.

Mathematical results concerning breather solutions
are based mostly on perturbation theory or variational methods,
except for special choices of the potentials $\phi$ and $\psi$
that admit simple solutions of a special form.
Surveys of rigorous results can be found in [\rAPankov,\rGalla].

Numerical methods have typically a much larger scope
and give more detailed information;
but they do not guarantee that the findings are correct, say up to small errors.
In this paper we give criteria that,
if satisfied by an approximate solution $\bar q$,
guarantee the existence of a true solution $q$ nearby.
For the solutions described in the theorem below,
the error $\|q-\bar q\|_\infty$ is shown to be less that $2^{-46}$.
Here, and in what follows,
$\|h\|_\infty=\sup_j\|h_j\|_\infty$ and $\|h_j\|_\infty=\sup_t|h_j(t)|$,
for any bounded function $(j,t)\mapsto h_j(t)$ on $\integer\times\real$.

\claim Theorem(breathers)
For each row in Table 1,
the equation \equ(Mainj) with the given parameter values
$(\phi,\psi,\omega)$ admits a $2\pi$-periodic solution $q$
that is real analytic in $t$, decreases at least exponentially in $|j|$,
and has norm $\|q\|_\infty>1$.
The entries $\pars$ and $\varpi$ in Table 1 describe symmetry properties.
The solution $q$ is symmetric or antisymmetric
with respect to time-reversal $t\mapsto-t$,
depending on whether $\pars=0$ or $\pars=1$, respectively.
An entry $\varpi=\pm$ indicates that $q$
admits a reflection symmetry as described by \equ(symmetry), with $\delta t=0$.

The entries $r$ and $\rho$ in Table 1 define bounds on the domain
of analyticity and decay rate of the solution; see Section 3.
The remaining entries will be described below.

Our proof of this theorem relies on estimates that are verified by a computer.
After writing \equ(Mainj) as a fixed point equation $G(q)=q$,
we prove \clm(breathers) by verifying that
a Newton-type map associated with $G$ is a contraction near $\bar q$.
This strategy has been used in many computer-assisted proofs, including [\rAKT,\rAKii].

\bigskip
{\eightpoint

\newcount\tmpnum 
\newdimen\tmpdim 
\def\opwarning#1{\immediate\write16{l.\the\inputlineno\space OPmac WARNING: #1.}}
\long\def\addto#1#2{\expandafter\def\expandafter#1\expandafter{#1#2}}
\long\def\isinlist#1#2#3{\begingroup \long\def\tmp##1#2##2\end{\def\tmp{##2}%
   \ifx\tmp\empty \endgroup \csname iffalse\expandafter\endcsname \else
                  \endgroup \csname iftrue\expandafter\endcsname \fi}
   \expandafter\tmp#1\endlistsep#2\end
}
\def\tabstrut{\strut}     
\def\tabiteml{\enspace}   
\def\tabitemr{\enspace}   
\def\vvkern{1pt}          
\def\hhkern{1pt}          


\newtoks\tabdata
\def\tabstrutA{\tabstrut}
\newcount\colnum
\def\ddlinedata{}
\def\vvleft{}

\def\table{\vbox\bgroup \catcode`\|=12 \tableA}
\def\tableA#1#2{\offinterlineskip \colnum=0 \def\tmpa{}\tabdata={}\scantabdata#1\relax
   \halign\expandafter{\the\tabdata\cr#2\crcr}\egroup}

\def\scantabdata#1{\let\next=\scantabdata
   \ifx\relax#1\let\next=\relax
   \else\ifx|#1\addtabvrule
      \else\isinlist{123456789}#1\iftrue \def\next{\scantabdataC#1}%
          \else \expandafter\ifx\csname tabdeclare#1\endcsname \relax
                \expandafter\ifx\csname paramtabdeclare#1\endcsname \relax
                   \opwarning{tab-declarator "#1" unknown, ignored}%
                \else \def\next{\expandafter \scantabdataB \csname paramtabdeclare#1\endcsname}\fi
             \else \def\next{\expandafter\scantabdataA \csname tabdeclare#1\endcsname}%
   \fi\fi\fi\fi \next
}
\def\scantabdataA#1{\addtabitem \expandafter\addtabdata\expandafter{#1\tabstrutA}\scantabdata}
\def\scantabdataB#1#2{\addtabitem\expandafter\addtabdata\expandafter{#1{#2}\tabstrutA}\scantabdata}
\def\scantabdataC {\def\tmpb{}\afterassignment\scantabdataD \tmpnum=}
\def\scantabdataD#1{\loop \ifnum\tmpnum>0 \advance\tmpnum by-1 \addto\tmpb{#1}\repeat
   \expandafter\scantabdata\tmpb
}

\def\unsskip{\ifdim\lastskip>0pt \unskip\fi}
\def\addtabitem{\ifnum\colnum>0 \addtabdata{&}\addto\ddlinedata{&\dditem}\fi
    \advance\colnum by1 \let\tmpa=\relax}
\def\addtabdata#1{\tabdata\expandafter{\the\tabdata#1}}
\def\addtabvrule{%
    \ifx\tmpa\vrule \addtabdata{\kern\vvkern}%
       \ifnum\colnum=0 \addto\vvleft{\vvitem}\else\addto\ddlinedata{\vvitem}\fi
    \else \ifnum\colnum=0 \addto\vvleft{\vvitemA}\else\addto\ddlinedata{\vvitemA}\fi\fi
    \let\tmpa=\vrule \addtabdata{\vrule}}

\def\crl{\crcr\noalign{\hrule}}
\def\crll{\crcr\noalign{\hrule\kern\hhkern\hrule}}

\def\crli{\crcr \omit
   \gdef\dditem{\omit\tablinefil}\gdef\vvitem{\kern\vvkern\vrule}\gdef\vvitemA{\vrule}%
   \vvleft\tablinefil\ddlinedata\crcr}
\def\crlli{\crli\noalign{\kern\hhkern}\crli}
\def\tablinefil{\leaders\hrule\hfil}

\def\crlp#1{\crcr \noalign{\kern-\drulewidth}%
   \omit \xdef\crlplist{#1}\xdef\crlplist{,\expandafter}\expandafter\crlpA\crlplist,\end,%
   \global\tmpnum=0 \gdef\dditem{\omit\crlpD}%
   \gdef\vvitem{\kern\vvkern\kern\drulewidth}\gdef\vvitemA{\kern\drulewidth}%
   \vvleft\crlpD\ddlinedata \global\tmpnum=0 \crcr}
\def\crlpA#1,{\ifx\end#1\else \crlpB#1-\end,\expandafter\crlpA\fi}
\def\crlpB#1#2-#3,{\ifx\end#3\xdef\crlplist{\crlplist#1#2,}\else\crlpC#1#2-#3,\fi}
\def\crlpC#1-#2-#3,{\tmpnum=#1\relax
   \loop \xdef\crlplist{\crlplist\the\tmpnum,}\ifnum\tmpnum<#2\advance\tmpnum by1 \repeat}
\def\crlpD{\global\advance\tmpnum by1
   \edef\tmpa{\noexpand\isinlist\noexpand\crlplist{,\the\tmpnum,}}%
   \tmpa\iftrue \kern-\drulewidth \tablinefil \kern-\drulewidth\else\hfil \fi}

\def\tskip{\afterassignment\tskipA \tmpdim}
\def\tskipA{\gdef\dditem{}\gdef\vvitem{}\gdef\vvitemA{}\gdef\tabstrutA{}%
    \vbox to\tmpdim{}\ddlinedata \crcr \noalign{\gdef\tabstrutA{\tabstrut}}}

\def\mspan{\omit \tabdata={\tabstrut}\let\tmpa=\relax \afterassignment\mspanA \mscount=}
\def\mspanA[#1]{\loop \ifnum\mscount>1 \csname span\endcsname \omit \advance\mscount by-1 \repeat
   \mspanB#1\relax}
\def\mspanB#1{\ifx\relax#1\def\tmpa{\def\tmpa####1}%
   \expandafter\tmpa\expandafter{\the\tabdata\ignorespaces}\expandafter\tmpa\else
   \ifx |#1\ifx\tmpa\vrule\addtabdata{\kern\vvkern}\fi \addtabdata{\vrule}\let\tmpa=\vrule
   \else \let\tmpa=\relax
      \ifx c#1\addtabdata{\tabiteml\hfil\ignorespaces##1\unsskip\hfil\tabitemr}\fi
      \ifx l#1\addtabdata{\tabiteml\ignorespaces##1\unsskip\hfil\tabitemr}\fi
      \ifx r#1\addtabdata{\tabiteml\hfil\ignorespaces##1\unsskip\tabitemr}\fi
   \fi \expandafter\mspanB \fi}

\newdimen\drulewidth  \drulewidth=0.4pt
\let\orihrule=\hrule  \let\orivrule=\vrule
\def\rulewidth{\afterassignment\rulewidthA \drulewidth}
\def\rulewidthA{\edef\hrule{\orihrule height\the\drulewidth}%
                \edef\vrule{\orivrule width\the\drulewidth}}

\long\def\frame#1{%
   \hbox{\vrule\vtop{\vbox{\hrule\kern\vvkern
      \hbox{\kern\hhkern\relax#1\kern\hhkern}%
   }\kern\vvkern\hrule}\vrule}}


\hfil\table{|c||c|c|c|c|c|c|c|c|c|c|c|c|c|}{\crl\tskip.1em
 label &$\phi_2$&$\phi_3$&$\phi_4$&$\psi_2$&$\psi_4$&$\omega$&$\varpi$&$\sigma$&$\pars$&$r$&$\rho$&$\ell$&stab\crli\tskip.1em
 \lab1  & $2^{-20}$ & $0$ & $2$     & $1/2$ & $0$   & $32/5$  &  $+$ & $0$ & $1$ & $2$     & $2$   &  $8$ & U \cr
 \lab2  & $2^{-20}$ & $0$ & $2$     & $1/2$ & $0$   & $32/5$  &  $-$ & $0$ & $1$ & $2$     & $2$   & $16$ & S \cr
 \lab3  & $1/8$    & $0$ & $1$     & $1/2$ & $0$   & $32/5$  &  $+$ & $1$ & $1$ & $2$     & $2$   & $15$ & S \cr
 \lab4  & $1/8$    & $0$ & $1$     & $1/2$ & $0$   & $32/5$  &  $-$ & $1$ & $1$ & $2$     & $2$   & $13$ & S \cr
 \lab5  & $-1$     & $0$ & $-27/16$& $1$   & $1$   & $32/5$  &  $-$ & $1$ & $1$ & $5/4$   & $5/4$ & $19$ & U \cr
 \lab6  & $-1$     & $0$ & $1$     & $0$   & $0$   & $256/5$ &  $-$ & $1$ & $1$ & $5/4$   & $5/4$ & $11$ & U \cr
 \lab7  & $-1$     & $0$ & $1$     & $0$   & $0$   & $256/5$ &  $+$ & $1$ & $1$ & $5/4$   & $5/4$ & $11$ & U \cr
 \lab8  & $-1/32$  & $0$ & $1$     & $1$   & $1/2$ & $32/5$  &  $-$ & $1$ & $1$ & $17/16$ & $2$   & $16$ & S \cr
 \lab9  & $-1/32$  & $0$ & $1$     & $1$   & $1/2$ & $32/5$  &  $+$ & $0$ & $1$ & $17/16$ & $2$   & $16$ & ns \cr
\lab10  & $-1/32$  & $0$ & $1$     & $1$   & $1/2$ & $32/5$  &  $+$ & $0$ & $1$ & $17/16$ & $2$   & $16$ & nu \cr
\lab11  & $-2^{-20}$& $0$ & $1$     & $1$   & $1/2$ & $32/5$  &  $+$ & $0$ & $1$ & $17/16$ & $2$   &  $6$ & U \cr
\lab12  & $1$      & $0$ & $1$     & $-1/2$& $0$   & $32/5$  &\snone& $1$ & $1$ & $9/8$   & $9/8$ & $21$ & U \cr
\lab13 & $1/4$    & $2$ & $4$     & $8$   & $0$   & $64/9$  &  $-$ & $1$ & $0$ & $9/8$   & $9/8$ & $15$ & S \crli}
}
\medskip
\centerline{\eightrm Table 1. parameter values and properties of solutions}
\bigskip

Our choice of parameters covers several different situations.
Attracting potentials $\phi$ and $\psi$
are used for the solutions $1$-$4$ and $13$.
The solutions $1$-$4$ cover the $4$ possible reflection symmetries
(bond-symmetry, bond-antisymmetry, site-symmetry, site-antisymmetry) with $\delta t=0$.
A globally repelling $\phi$ is used for solution $5$.
The solutions $6$ and $7$ correspond to typical mountain pass configuration
for the functional $\Lagrangian$ described below.
For the solutions $8$-$11$ we have a coercive potential $\phi$,
which is repelling in a neighborhood of $0$, while $\psi$ is attracting.
The opposite is the case for solution $12$.
Here we chose a solution with four bumps with different fundamental time-periods.
Finally, solution $13$ is an even function of time
(while our other solutions are odd) and has a nonzero time-average.

The solutions $11$ and $12$ are shown in Fig.~1,
at sites $j\in\integer$ where $\|q_j\|_\infty>2^{-10}\|q\|_\infty$.
The values for non-integer $j$ are obtained by linear interpolation.
Graphs for all solutions,
  in the form of animations, can be found in [\rFiles].

\demo Remark(negbreather)
If $\phi$ and $\psi$ both even,
then $-q$ solves the equation \equ(Mainj) whenever $q$ does.
This applies to our solutions $1$-$12$.

\smallskip
In order to discuss the stability of these solutions,
we write the second order equation \equ(Mainj) for $q$
as a first order equation for the pair
$\buu=[{q\atop p}\bigr]$, where $p=\dot q$.
The resulting equation is in fact Hamiltonian,
with the Hamiltonian function given by
$$
\Hamiltonian(q,p)=\sum_j\thalf(p_j)^2+\Potential(q)\,,\quad
\Potential(q)=\omega^{-2}\sum_j\Bigl[\phi(q_{j+1}-q_j)+\psi(q_j)\Bigr]\,.
\equation(Hamiltonian)
$$
In other words, $\dot q_j=\partial_{p_j}\Hamiltonian$
and $\dot p_j=-\partial_{q_j}\Hamiltonian$.
The corresponding time-$t$ map $\buu(0)\mapsto\buu(t)$ will be denoted by $\FLOW_t$.
Since we are interested in breather solutions,
it suffices to consider initial conditions $\buu(0)$
whose components $q(0)$ and $p(0)$ belong to $\HH=\ell^2(\integer)$.
Let now $\buu$ be a $2\pi$-periodic orbit in $\HH^2$.
Then each $\FLOW_t$ is well-defined and differentiable
in some open neighborhood of $\buu(0)$ in $\HH^2$.
Define $\Flow(t)$ to be the derivative $D\Psi_t(\buu(0))$.
We say that the orbit $\buu$ is spectrally stable
if the spectrum of $\Flow(2\pi)$ belongs to the closed unit disk.
In fact, we can replace ``unit disk'' by ``unit circle'',
since the spectrum is invariant under $z\mapsto\bar z$ and $z\mapsto z^{-1}$,
due to the Hamiltonian nature of the flow.

The spectrum of the time-$2\pi$ map for trivial solution $\buu=0$
is easily seen to be the set of all complex numbers $e^{2\pi iz}$
for which $z^2$ is real and belongs to
the interval bounded by $\omega^{-2}\psi_2$ and $\omega^{-2}[\psi_2+4\phi_2]$.
It is not hard to see that this set $\Sigma^e$
also constitutes the essential spectrum
of the time-$2\pi$ map $\Flow(2\pi)$ for our breather solution.

\demo Remark(translation)
If $\buu$ is an orbit for the flow generated by $\Hamiltonian$,
then so is $t\mapsto\buu(t-c)$, for any constant $c$.
This implies e.g.~that $\dot\buu$ is an eigenvector of $\Flow(2\pi)$
with eigenvalue $1$.

\claim Theorem(stability)
Consider the solution described in \clm(breathers),
associated with one of the rows of Table 1.
If the entry in the last column is a {\rm``S''} or {\rm``U''},
then the solution is spectrally stable or unstable, respectively.

An entry ``ns'' or ``nu'' in Table 1
means that the solution appears to be spectrally stable or unstable, respectively,
based on numerical results.
We did not succeed in validating these results, due to a limited ability
of our methods to deal with continuous spectrum.
In particular, solution 9 appears to have an eigenvalue with the ``wrong''
Krein signature embedded in the continuous spectrum.
This is a notoriously difficult situation.

The first step in our proof of \clm(stability)
is to show that it suffices to work with truncated systems
whose time-$2\pi$ maps $\Flow_n(2\pi)$ are essentially matrices.
This result should be of independent interest.
The spectrum of $\Flow(2\pi)$
outside $\Sigma^e$ consist of isolated eigenvalues with finite multiplicities.
In Theorem 2.13, we show that these eigenvalues are approximated
by the eigenvalues of $\Flow_n(2\pi)$ for large $n$.

For the part of \clm(stability) that deals with spectral stability,
we use ideas from [\rAKii], where spectral (in)stability was proved
for some solutions of a periodically perturbed wave equation.
The proof of instability is complicated by the presence
of the above-mentioned eigenvalue $1$.
As a consequence, our instability results are restricted to cases
where the continuous spectrum of $\Flow(2\pi)$ is either very narrow
(solutions $1$ and $11$)
or includes a real interval (solutions $5$, $6$, $7$, and $12$).
In the first case, linear instability is the result of eigenvalues outside the unit disk.
By standard results on invariant manifolds,
this implies e.g.~that the solutions $1$ and $11$ are truly (not just linearly) unstable.

Solution $3$ and its spectrum are shown in Fig.~2.
Spectrum that is not marked with dots lies in the pink and gray arcs. The color
indicates the Krein signature: red or pink means positive, black or gray means negative.
Blue (cyan) dashes mark the primary (non-primary) separating values; see Section 4.
Spectral plots for some of our other solutions can be found in [\rFiles].

\medskip
For each of the solutions $q$ described in \clm(breathers),
there exists a finite lattice interval
$\JJ=\{j\in\integer: \sigma-\ell\le j-J\le\ell\}$
which we call the ``approximate support'' of $q$.
The value of $\ell$ is given in Table 1.

\claim Theorem(multi)
Consider a fixed choice of parameters $(\phi,\psi,\omega,\sigma,\pars,r,\rho)$
from Table 1, excluding rows 5 and 13.
Let $m\mapsto q^m$ be a sequence (finite or infinite)
of solutions of the equation \equ(Mainj),
associated with these parameter values,
as described in \clm(breathers) and \dem(negbreather).
By considering translates,
we assume now that the approximate supports $\JJ_m$ of these solutions
are mutually disjoint.
Assume in addition that the distance between
any two adjacent approximate supports is even.
Then there exists a solution $q$ of \equ(Mainj),
with the property that $\|q_j-q^m_j\|_\infty<2^{-45}$
whenever $j\in\JJ_m$ for some $m$,
and $\|q_j\|_\infty<2^{-50}$ whenever $j\not\in\JJ_m$ for all $m$.

The idea of the proof is of course to use that
the breathers $q^m$ interact very little if they are placed sufficiently far apart.
This idea has been used e.g.~in [\rAGT,\rFLMii,\rPS]
to construct and analyze extended solutions for the FPU model and other lattice systems.
Notice however that the notion of ``sufficiently far''
in \clm(multi) is specific and very mild.

Our proof of this theorem is based again on a contraction mapping argument.
Here we have to work with $\ell^\infty$ type spaces;
but since the interactions have finite range,
the estimates that are needed are not much stronger
than what is required for our proof of \clm(breathers).
Nevertheless, these estimates fail for the solutions $5$ and $13$.
Instead of treating these cases differently,
we chose to exclude them from \clm(multi), just for simplicity.

We expect that a multi-breather $q$ is spectrally stable only in very special cases.
Assuming that the sets $\JJ_m$ are placed sufficiently far apart,
each of the breathers $q^m$ will have to be spectrally stable.
But this is not sufficient, since a system of $N$ non-interacting breathers
has an eigenvalue $1$ with multiplicity $2N$.
Under the influence of a small interaction,
all but $2$ of these eigenvalues can move away from $1$.
It may be possible to keep these eigenvalues on the unit circle
by using time-translates $q^m(\bdot-c_m)$ in the construction of $q$,
with properly chosen constants $c_m$.
But this is outside the scope of our current methods.

As mentioned earlier,
FPU-type models are accessible to variational methods as well.
Time-periodic solutions of \equ(Mainj) with period $2\pi$
can be found as critical points of the Lagrangian functional
$$
\Lagrangian(q)=\int_0^{2\pi}\biggl[\,\sum_j\thalf(\dot q_j)^2
-\Potential(q)\biggr]\,dt\,,
\equation(Lagrangian)
$$
defined on a suitable Hilbert space of functions $q:\rmS^1\to\ell^2$.
We refer to [\rAGb,\rAGa] for early results and to [\rAPankov] for a survey.
The coefficients $\phi_i$ and $\psi_i$ determine the geometry of the functional.
In particular, the choice $\psi=0$ and $\phi_2<0<\phi_4$
yields a mountain pass geometry.
This is the first (and simplest) case that was considered with critical point theory,
and it is the only case for which the existence of multi-breather solutions
has been proved variationally [\rAGT].

\medskip
The remaining part of this paper is organized as follows.
In Section 2 we consider the time-$2\pi$ map
for an infinite chain and its spectral approximation
by time-$2\pi$ maps $\Flow_n(2\pi)$ for chains of length $2n$.
Section 3 is devoted to the task of
proving Theorems \clmno(breathers) and \clmno(multi).
Our proof of these theorems requires estimates on approximate solutions.
The same is true for our proof of \clm(stability),
which is given in Section 4.
The instability proof (for solutions $1$ and $11$) uses a perturbation argument.
The stability proof (for solutions $2$, $3$, $4$, $8$, and $13$) uses Krein signatures,
and a monotonicity argument from [\rAKii],
to control the eigenvalues of $\Flow_n(2\pi)$.
The estimates that are needed in Sections 3 and 4 are proved with the aid of a computer;
a rough description is given in Section 5,
and for details we refer to the source code of our programs [\rFiles].

\vskip1.0cm
\hbox{\hskip0.0cm
\includegraphics[height=4.0cm,width=5.5cm]{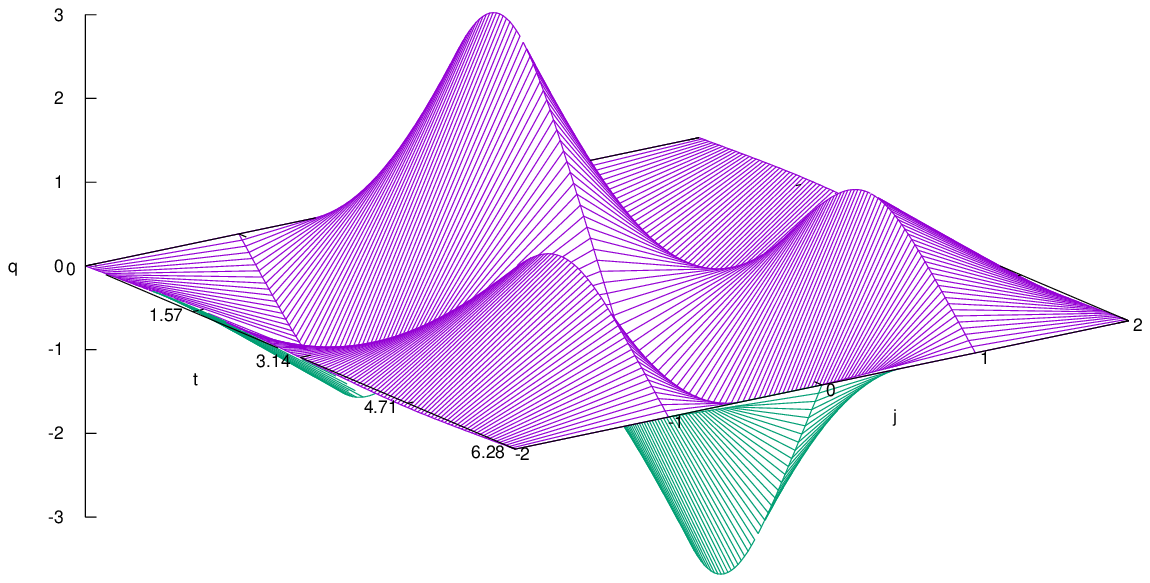}
\includegraphics[height=4.0cm,width=8.0cm]{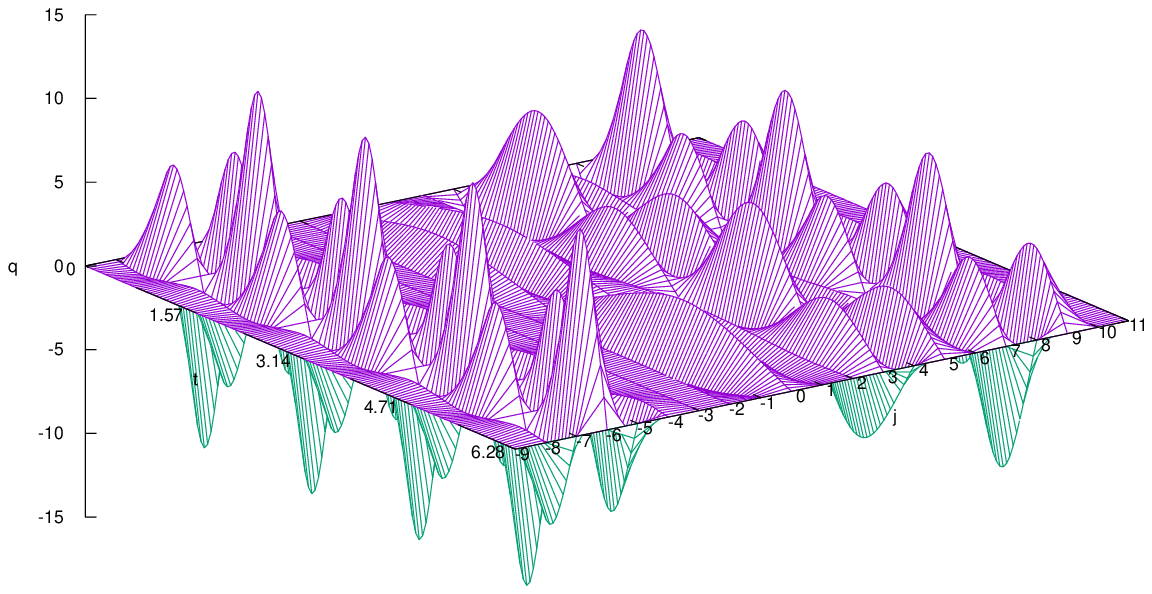}}
\vskip0.4cm
\centerline{\eightpoint{\noindent\bf Figure 1.}
Solutions $11$ and $12$.}

\vskip1.0cm
\hbox{\hskip0.0cm
\raise0.0cm\hbox{\includegraphics[height=5.0cm,width=7.0cm]{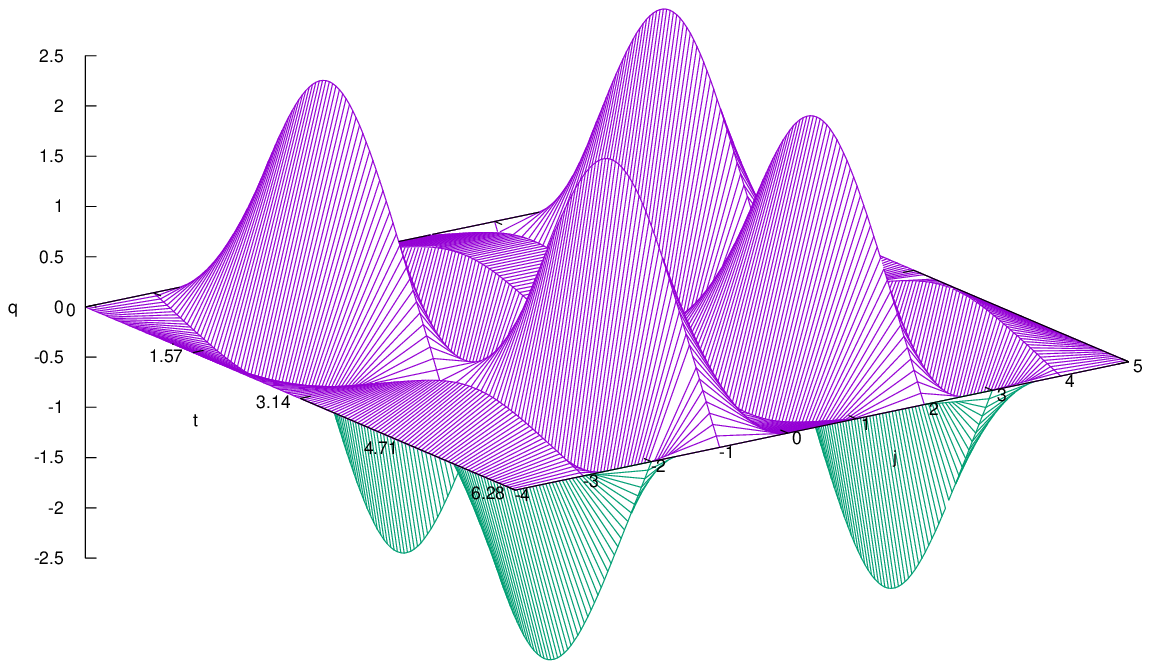}}
\hskip0.5cm
\includegraphics[height=5.5cm,width=5.5cm]{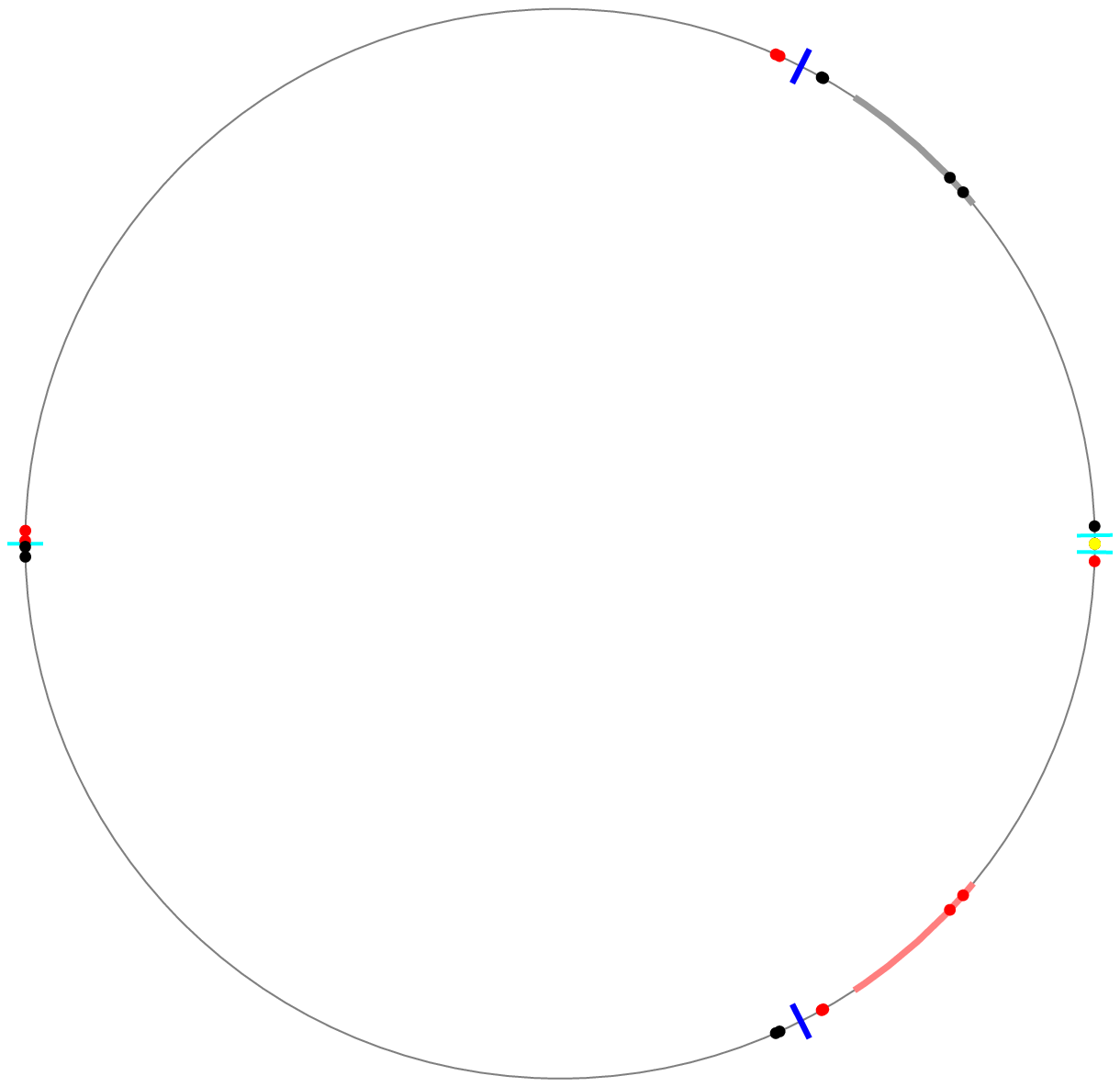}}
\vskip0.4cm
\centerline{\eightpoint{\noindent\bf Figure 2.}
Solution $3$ and its spectrum.}

\section Infinite chains and approximations

After introducing some notation,
we discuss spectral properties of operators
in a class that includes the linearized time-$2\pi$ maps
for exponentially decreasing time-periodic breathers.

\subsection Notation

By a chain $q$ we mean a real-valued function $j\mapsto q_j$ on $\integer$.
Complex-valued functions will be considered only for spectral theory.
If $f$ is any real-valued function on $\real$,
then $f(q)$ denotes the chain with values $f(q)_j=f(q_j)$.
We say that a chain $q$ is site-centered (bond-centered)
if $q$ has a symmetry \equ(symmetry) with $\sigma=0$ ($\sigma=1$).
This restriction to $\sigma=0$ and $\sigma=1$ is motivated mainly
by computational simplicity.

Given $\sigma\in\{0,1\}$, we set
$$
(\nablas q)_j=q_{j+\sigma}-q_{j+\sigma-1}\,,
\equation(nablasigmaDef)
$$
for every chain $q$ and every integer $j$.
This defines a continuous linear operator $\nablas$
on $\HH=\ell^2(\integer)$.
Its adjoint is given by $\nablas^\ast=-\nabla_{\!1-\sigma}$.
The equation \equ(Mainj) can now be written as
$$
\ddot q=-K(q)\,,\qquad
K(q)=\omega^{-2}\bigl[\nablas^\ast\phi'(\nablas q)+\psi'(q)\bigr]\,.
\equation(Main)
$$
The corresponding first-order equation is
$$
\dot\buu=Y(\buu)\,,\qquad\buu=\stwovec{q}{p}\,,\qquad
Y(\buu)=\stwomat{0}{\id}{-K(q)}{0}\,.
\equation(uFlow)
$$

We note that the operator $K$ is independent of the choice of $\sigma$.
Our reason for considering two distinct versions of the lattice gradient is
that $\nablao$ maps site-centered chains to bond-centered chains,
while $\nablai$ maps bond-centered chains to site-centered chains.

Let now $\buu$ be a fixed $2\pi$-periodic orbit for the flow \equ(uFlow),
and define
$$
\alpha=\omega^{-2}\phi''(\nablas q)\,,\qquad\beta=\omega^{-2}\psi''(q)\,.
\equation(abDef)
$$
Consider an orbit $\buu+\buv$ close to $\buu$.
To first order in $\buv$, we have
$$
\dot\buv=X(\buu)\buv\,,\qquad
X(\buu)\defeq DY(\buu)=\stwomat{0}{\id}{-H(q)}{0}\,,
\equation(vFlow)
$$
where $H(q)$ is the linear operator
$$
H(q)v\defeq DK(q)v=\nablas^\ast\alpha\nablas v+\beta v\,,\qquad v\in\HH\,.
\equation(HqDef)
$$
If $\FLOW_t$ denotes the time-$t$ map for the flow \equ(uFlow),
then the time-$t$ map for the flow \equ(vFlow)
is given by the linear operator $\Flow(t)=D\FLOW_t(\buu(0))$.
Our goal is to analyze the spectrum of the operator $\Flow(2\pi)$.

Notice that $H(q(t))$ is $2\pi$-periodic in $t$ and self-adjoint for each $t$.
Furthermore, far from the origin, where $q$ is close to zero,
$H(q)$ is well approximated by the operator
$$
H(0)=\bar\alpha\nablas^\ast\nablas+\bar\beta\,.
\equation(HoDef)
$$
Here, $\bar\alpha=\omega^{-2}\phi_2$ and $\bar\beta=\omega^{-2}\psi_2$.

In what follows, we consider operators of the type \equ(HqDef) and \equ(HoDef),
where $\alpha$ and $\beta$ can be more general $2\pi$-periodic curves in $\HH$.
But we will assume that $\tilde\alpha=\alpha-\bar\alpha$
and $\tilde\beta=\beta-\bar\beta$ decrease exponentially.
To simplify notation, the remaining part of this section
is formulated for $\nablao$ only.

\subsection Some spaces

In order to discuss the spectrum of $\Flow(2\pi)$,
we will need certain constructions for several different Hilbert spaces.
To avoid undue repetition,
we start by considering a fixed but arbitrary Hilbert space $\buh$.
The inner product in $\buh$ will be denoted by $\langle\bdot\,,\bdot\rangle$.
Here, and in what follows, an inner product
is always assumed to be linear in its second argument,
and antilinear in the first argument.
Consider the vector space $\buc_0(\integer,\buh)$ of all sequences
$v:\integer\to\buh$ with the property that $v_j=0$
for all but finitely many $j\in\integer$.
For every $\rho\ge 0$, we define $\buH_\rho(\buh)$ to be the Hilbert space
obtained as the closure of $\buc_0(\integer,\buh)$ with respect to the norm
$$
\|v\|_\rho={\textstyle\sqrt{\langle v,v\rangle_\rho}}
\,,\qquad\langle v,v'\rangle_\rho
=\sum_{j\in\integer}\cosh(2j\rho)\langle v_j,v'_j\rangle\,.
\equation(HHrhoIprod)
$$
The equation
$$
\bigl(U(\rho)v\bigr)_j=\cosh(2\rho j)^{1/2}v_j\,,
\qquad v\in\buH_0(\buh)\,,\quad j\in\integer\,,
\equation(UirhoDef)
$$
defines a unitary operator $U(\rho)$ from $\buH_\rho(\buh)$ to $\buH_0(\buh)$.
The Fourier transform $\FF v$ of a function $v\in\buH_0(\buh)$ is given by
$$
\FF v=\tilde v\,,\qquad
\tilde v(\varphi)=\sum_{j\in\integer} v_je^{-i\varphi j}\,.
\equation(FourierSeries)
$$
The Fourier transform $\FF$ defines a unitary operator
from $\buH_0(\buh)=\ell^2(\integer,\buh)$ to the Hilbert space $\rmL^2(\II,\buh)$
with the inner product
$$
\langle g,f\rangle_{\sss\rmL^2}
={1\over|\II|}\int_\II\Langle g(\varphi),f(\varphi)\Rangle\,d\varphi\,,
\qquad \II=[-\pi,\pi]\,.
\equation(tildeHHIprod)
$$
Consider now $v\in\buH_\rho(\buh)$ with $\rho>0$.
Then the sum in \equ(FourierSeries) extends $\tilde v$
to an analytic function on the interior of the strip
$$
S_\rho=\bigl\{\varphi\in\complex: |\Im\varphi|\le\rho\bigr\}\,.
\equation(rhoStrip)
$$
Furthermore, the function $\varphi\mapsto\tilde v(\varphi+i\rho)$
belongs to $\rmL^2(\II,\buh)$, and
$$
\|v\|_\rho^2=\thalf\|\tilde v(\bdot+i\rho)\|_{\sss\rmL^2}^2
+\thalf\|\tilde v(\bdot-i\rho)\|_{\sss\rmL^2}^2\,.
\equation(HHrrhoNorm)
$$
We will also need to approximate functions in $\buH_\rho(\buh)$
by functions that are supported in $\integer_n=\{j\in\integer:-n<j\le n\}$
for some positive integer $n$.
To this end, define a projection $P_n:\buH_0(\buh)\to\buH_0(\buh)$ by setting
$$
(P_n v)_j=\cases{v_j &if $j\in\integer_n$,\cr 0 &otherwise,\cr}
\qquad v\in\buH_0(\buh)\,,\quad j\in\integer\,.
\equation(PnDef)
$$
The following proposition is an immediate consequence
of the fact that, if $0\le\rho\le\varrho$, then
$\cosh(2\varrho n)\|(\id-P_n)v\|_\rho\le\cosh(2\rho n)\|(\id-P_n)v\|_\varrho$.

\claim Proposition(TailBound)
Let $0\le\rho<\varrho$.
Assume that $v\in\buH_\varrho(\buh)$ is nonzero
and satisfies a bound $\|v\|_\varrho\le C\|v\|_\rho$.
If $n$ is sufficiently large, so that
$\cosh(2\varrho n)>C^2\cosh(2\rho n)$, then
$$
\|(\id-P_n)v\|_\rho
\le\biggl({\cosh(2\rho n)(C^2-1)\over
\cosh(2\varrho n)-C^2\cosh(2\rho n)}\biggr)^{1/2}
\|P_nv\|_\rho\,.
\equation(TailBound)
$$

\medskip
The subspace $P_n\buH_\rho(\buh)$ of $\buH_\rho(\buh)$
can be identified with the Hilbert space $\buH_{\rho,n}(\buh)$ of all functions
$v:\integer_n\to\buh$, equipped with the inner product
$$
\langle v,v'\rangle_\rho=\sum_{j\in\integer_n}\cosh(2\rho j)\langle v_j,v'_j\rangle\,.
\equation(HHnrhoIprod)
$$
Notice that $\buH_{\rho,n}(\buh)$ agrees with $\buH_{0,n}(\buh)$
as a vector space, for any $\rho>0$.
The Fourier transform $\tilde v=\FF_n v$
of a function $v\in\buH_{0,n}(\buh)$ is defined as in \equ(FourierSeries),
but with the sum ranging over $j\in\integer_n$ only.
The Fourier transform $\FF_n$ defines a unitary operator
from $\buH_{0,n}(\buh)=\ell^2(\integer_n,\buh)$ to the Hilbert space $\ell^2(\II_n,\buh)$
with the inner product
$$
\langle g,f\rangle_{\ell^2}
={1\over|\II_n|}\sum_{\varphi\in\II_n}\Langle g(\varphi),f(\varphi)\Rangle\,,
\qquad\II_n={\pi\over 2n}\integer_n\,.
\equation(tildeHHnIprod)
$$

\claim Notation(OpNorm)
The operator norm of a continuous linear operator
$L$ on $\buH_\rho(\buh)$ or $\buH_{\rho,n}(\buh)$ will be denoted by $\|L\|_\rho$.

\subsection The flow and its spectrum

Let $\HH_\rho=\buH_\rho(\complex)$,
using the inner product $\langle z,z'\rangle=\bar z z'$ on $\complex$.
Here, and in what follows, $\rho$ is a fixed but arbitrary
nonnegative real number, unless specified otherwise.

Consider the space $\HH_\rho^2$ of all pairs $\buv=\bigl[{v\atop\nu}\bigr]$
with components $v,\nu\in\HH_\rho$.
We will identify $\HH_\rho^2$ with the space $\buH_\rho\bigl(\complex^2\bigr)$
of sequences $j\mapsto\buv_j=\bigl[{v_j\atop\nu_j}\bigr]$,
using the inner product $\langle\buv_j,\buv_j'\rangle=\bar v_j v'_j+\bar\nu_j\nu'_j$
on $\complex^2$. On $\HH_\rho^2$ we consider the operator
$$
X_0=\stwomat{0}{\id}{-H_0}{0}\,,\qquad
H_0=\bar\alpha\nablao^\ast\nablao+\bar\beta\,,
\equation(XoDef)
$$
where $\bar\alpha$ and $\bar\beta$ are fixed but arbitrary real numbers.
Notice that the operator $\FF H_0\FF^{-1}$ on $\FF\HH_\rho$
is multiplication by the function $h_0$,
$$
h_0(\varphi)=2\bar\alpha\bigl[1-\cos(\varphi)\bigr]+\bar\beta\,.
\equation(hoDef)
$$
So the spectrum of the operators $H_0:\HH_\rho\to\HH_\rho$
and $X_0:\HH_\rho^2\to\HH_\rho^2$ are given by
$R_\rho$ and $iR_\rho^{1/2}$, respectively, where
$$
R_\rho=\range\bigl(h_0:S_\rho\to\complex\bigr)\,,\qquad
iR_\rho^{1/2}=\bigl\{z\in\complex: -z^2\in R_\rho\bigr\}\,,
\equation(SpecHoHHrho)
$$
and where $S_\rho$ is the strip defined in \equ(rhoStrip).

A straightforward computation shows that
$$
e^{tX_0}=\twomat{\cos(t\bmy)}{\bmy^{-1}\sin(t\bmy)}
{-\bmy\sin(t\bmy)}{\cos(t\bmy)}\,,
\qquad\bmy=H_0^{1/2}\,.
\equation(Phioy)
$$
The choice of square root does not matter as this point,
since $\cos(t\bmy)$ and $\bmy^{\pm 1}\sin(t\bmy)$
are even functions of $\bmy$.

\bskip
Let $(j,t)\mapsto\tilde\alpha_j(t)$ and $(j,t)\mapsto\tilde\beta_j(t)$
be two functions on $\integer\times\real$ with the following properties.

\claim Condition(alphabetaProp)
The functions $\tilde\alpha_j$ and $\tilde\beta_j$ are
continuous and $2\pi$-periodic.
Furthermore, the sequences $U(\rho_0)\tilde\alpha(t)$
and $U(\rho_0)\tilde\beta(t)$ are bounded for some $\rho_0>0$,
uniformly in $t$.
Here $U(\rho_0)$ denotes the multiplication operator
as defined in \equ(UirhoDef).

This defines two families of compact linear operators
on $\HH_\rho$ via pointwise multiplication:
$\bigl(\tilde\alpha(t)v\bigr)_j=\tilde\alpha_j(t)v_j$ and
$\bigl(\tilde\beta(t)v\bigr)_j=\tilde\beta_j(t)v_j$.
Define also
$$
H_1(t)v=\nablao^\ast\tilde\alpha(t)\nablao v+\tilde\beta(t)v\,,\qquad t\in\real\,.
\equation(HiDef)
$$
On $\HH_\rho^2$ we consider the flow given by the equation
$$
{d\over dt}\buv(t)=X(t)\buv(t)\,,\qquad
X(t)=\stwomat{0}{\id}{-H(t)}{0}\,,\quad
H(t)=H_0+H_1(t)\,.
\equation(XFlow)
$$
The corresponding time-$t$ map $\buv(0)\mapsto\buv(t)$ will be denoted by $\Flow(t)$.
Clearly, $\Flow(t)$ is bounded on $\HH_\rho^2$ for each $t\in\real$.
In order to get more detailed information on $\Flow(t)$,
we use that $\Flow$ can be obtained by solving the Duhamel equation
$$
\Flow(t)=e^{tX_0}
+\int_0^t e^{(t-s)X_0}\bigl[X(s)-X_0\bigr]\Flow(s)\,ds\,.
\equation(PhiIntEqu)
$$
Notice that $X(s)-X_0=P_2^\ast H_1(s)P_1$,
where $P_1=[1\ 0]$ and $P_2=[0\ 1]$ are the operators
that assign to a vector $\buv=\bigl[{v\atop\nu}\bigr]$
in $\HH_\rho^2$ its component $v\in\HH_\rho$ and $\nu\in\HH_\rho$, respectively.
Since $H_1$ is a continuous curve of compact linear operators,
the integral in this equation defines a compact linear operator.
Thus, $\Flow(t)$ is a compact perturbation of $e^{tX_0}$.
In particular, the essential spectrum of $\Flow(2\pi)$
agrees with the essential spectrum of $e^{2\pi X_0}$
which is
$$
\Sigma_\rho^e=\exp\bigl(\,2\pi iR_\rho^{1/2}\,\bigr)\,.
\equation(EssSpec)
$$
Notice that $\Sigma_0^e$ is included in the union
of the unit circle and the real line.
If $\Sigma_0^e$ does not cover the entire circle,
then the complement of $\Sigma_0^e$ is connected.
The same holds for $\Sigma_\rho^e$,
if $\rho>0$ is chosen sufficiently small.

\claim Theorem(PhiEigen)
Let $\rho\ge 0$, and assume that the complement of $\Sigma_\rho^e$
is connected.
Then the spectrum of $\Flow(2\pi):\HH_\rho^2\to\HH_\rho^2$
outside $\Sigma_\rho^e$ consists of eigenvalues
with finite (algebraic) multiplicities.
These eigenvalues can accumulate only at $\Sigma_\rho^e$.

This theorem is a consequence of the following fact.
A bounded linear operator with essential spectrum $\Sigma^e$
has at most countable spectrum in the unbounded component of $\complex\setminus\Sigma^e$.
The spectrum in this component consists of isolated eigenvalues
with finite (algebraic) multiplicities.
See e.g.~the proof of Theorem 4.5.33 in [\rKato].

\subsection Eigenvalues and eigenvectors

Multiplying both sides of the equation \equ(PhiIntEqu) from the left by $e^{-tX_0}$,
we find that the family of operators $A(t)=e^{-tX_0}\Flow(t)-\id$
satisfies the equation
$$
A(t)=-\int_0^t B(s)\bigl[\id+A(s)\bigr]\,ds\,,
\equation(AIntEqu)
$$
where
$$
\eqalign{
B(s)
&=e^{-sX_0}P_2^\ast H_1(s)P_1 e^{sX_0}\cr
&=\twomat{-\bmy^{-1}\sin(s\bmy)H_1(s)\cos(s\bmy)}
{-\bmy^{-1}\sin(s\bmy)H_1(s)\bmy^{-1}\sin(s\bmy)}
{\cos(s\bmy)H_1(s)\cos(s\bmy)}
{\cos(s\bmy)H_1(s)\bmy^{-1}\sin(s\bmy)}\,.\cr}
\equation(BsDef)
$$
The integral equation \equ(AIntEqu) for $A$ can be solved by iteration,
$$
A(t_0)=\sum_{n=1}^\infty(-1)^n\int_0^{t_0}\!\!dt_1\int_0^{t_1}\!\!dt_2\cdots\int_0^{t_{n-1}}\!dt_n
B(t_1)B(t_2)\cdots B(t_n)\,.
\equation(AtoIter)
$$

Denote by $\rho_0$ the decay rate of the sequences
$\tilde\alpha(t)$ and $\tilde\beta(t)$, as described in \clm(alphabetaProp).
Then the following holds.

\claim Proposition(CompactOps)
Assume that $0<\rho<\rho_0$.
Then $H_1(t)$ defines a compact linear operator from $\HH_0$ to $\HH_\rho$.
Furthermore, $B(t)$ and $A(t)$ define
compact linear operators from $\HH_0^2$ to $\HH_\rho^2$.
The bounds are uniform in $t$, if $t$ is restricted to a bounded interval.

\proof
Using that $U(\rho)^{-1}\nabla U(\rho)$ is bounded on $\HH_\rho$,
and that $\tilde\alpha(t)U(\rho)$ and $\tilde\beta(t)U(\rho)$
are compact on $\HH_\rho$, we see that
$$
H_1(t)U(\rho)=\nabla^\ast\tilde\alpha(t)U(\rho)\bigl[U(\rho)^{-1}\nabla U(\rho)\bigr]
+\tilde\beta(t)U(\rho)
\equation(HitUirho)
$$
defines a compact operator on $\HH_\rho$.
This proves the first claim.

Composing $\cos(t\bmy):\HH_0\to\HH_0$ with $H_1:\HH_0\to\HH_\rho$
and $\cos(t\bmy):\HH_\rho\to\HH_\rho$,
we see that
$$
B_{1,1}(t)=\cos(t\bmy)H_1(t)\cos(t\bmy)
\equation(Biit)
$$
defines a compact linear operator from $\HH_0$ to $\HH_\rho$.
Similarly for the other components of the operator $B(t)$ defined in \equ(BsDef).
Thus, $B(t)$ is compact as a linear operator from $\HH_0^2$ to $\HH_\rho^2$.
By \equ(AtoIter) the same is true for $A(t)$.
\qed

\claim Corollary(EigenDecay)
Let $0<\rho<\rho_0$.
Let $\buv\in\HH_0^2$ be an eigenvector of $\Flow(2\pi)$
with eigenvalue $\lambda\in\complex\setminus\Sigma_\rho^e$.
Then $\buv$ belongs to $\HH_\rho^2$, and
$$
\buv=\bigl[\lambda\id-e^{2\pi X_0}\bigr]^{-1}e^{2\pi X_0}A(2\pi)\buv\,.
\equation(DecayOne)
$$

\proof
By assumption we have
$\lambda\buv=\Flow(2\pi)\buv=e^{2\pi X_0}\bigl[\id+A(2\pi)\bigr]\buv$,
and thus
$$
\bigl[\lambda\id-e^{2\pi X_0}\bigr]\buv
=e^{2\pi X_0}A(2\pi)\buv\,.
\equation(DecayTwo)
$$
By \clm(CompactOps), the right hand side of this equation
belongs to $\HH_\rho^2$.
And $\lambda\id-e^{2\pi X_0}$ has a bounded inverse on $\HH_\rho^2$
since $\lambda\not\in\Sigma_\rho^e$. This proves the claim.
\qed

For completeness, we state a partial generalization of \clm(EigenDecay),
which follows from (the proof of) Proposition 4.3 in [\rAKii].

\claim Proposition(SpecDecay)
Let $\lambda$ be an isolated eigenvalue of $\Flow(2\pi):\HH_0^2\to\HH_0^2$.
Then each vector in the corresponding spectral subspace
belongs to $\HH_\rho^2$ for some $\rho>0$.

\subsection Resolvent estimates

Consider now the case where $\bar\beta$ and $\bar\beta+4\bar\alpha$
are both nonnegative.
Define $\bmy=H_0^{1/2}$ by using a square root
function that is analytic in $\complex\setminus(-\infty,0]$.
Then a partial diagonalization of $X_0$ is given by
$$
X_0=\twomat{0}{\id}{-\bmy^2}{0}
=\Lambda\twomat{-i\bmy}{0}{0}{i\bmy}\Lambda^{-1}\,,
\equation(XoPartialDiag)
$$
where
$$
\Lambda
={1\over\sqrt{2}}\twomat{\bmy^{-1/2}}{\bmy^{-1/2}}{-i\bmy^{1/2}}{i\bmy^{1/2}}\,,
\qquad
\Lambda^{-1}
={1\over\sqrt{2}}\twomat{\bmy^{1/2}}{i\bmy^{-1/2}}{\bmy^{1/2}}{-i\bmy^{-1/2}}\,.
\equation(XoPartialDiagLambda)
$$

\claim Definition(rhoo)
Let $\rho_1$ be the largest positive real number
with the property that the range of $h_0:S_\rho\to\complex$ is contained
in $\complex\setminus(-\infty,0]$, whenever $\rho<\rho_1$.
For simplicity, we assume from now on that $\rho_0\le\rho_1$.

Clearly, if $\rho<\rho_1$,
then $\Lambda$ and $\Lambda^{-1}$ define bounded linear operators on $\HH_\rho^2$.

\claim Proposition(ResBound)
Let $0\le\rho<\rho_1$ and define $C=\|\Lambda\|_\rho\|\Lambda^{-1}\|_\rho$.
Then
$$
C\bigl\|e^{2\pi X_0}\buv-z\buv\bigr\|_\rho
\ge\,\dist\bigl(z,\Sigma_\rho^e\bigr)\|\buv\|_\rho\,,
\qquad\buv\in\HH_\rho^2\,,\quad z\in\complex\,.
\equation(ResBound)
$$

\proof
Let $u\in\HH_\rho$,
and consider the Fourier transform $\tilde u$
for arguments $\varphi+i\kappa$ with $\varphi,\kappa\in\real$
and $|\kappa|\le\rho$.
On this domain we have the bound
$$
\Bigl|\Bigl(e^{\pm 2\pi ih_0^{1/2}}-z\Bigr)\tilde u\Bigr|
\ge\dist\bigl(z,\Sigma_\rho^e\bigr)|\tilde u|\,,
\equation(ResBoundFive)
$$
for any choice of the square root function.
Using the identity \equ(HHrrhoNorm),
this implies the bound
$$
\bigl\|\bigl(e^{\pm 2\pi i\sbmy}-z\id\bigr)u\bigr\|_\rho
\ge\dist\bigl(z,\Sigma_\rho^e\bigr)\|u\|_\rho\,.
\equation(ResBoundSeven)
$$
By \equ(XoPartialDiag) we have
$$
\twomat{e^{-2\pi i\sbmy}-z\id}{0}{0}{e^{2\pi i\sbmy}-z\id}
=\Lambda^{-1}\bigl[e^{2\pi X_0}-z\id\bigr]\Lambda\,.
\equation(ResBoundOne)
$$
Let $\buv\in\BB_\rho$ and $\buv'=\Lambda^{-1}\buv$. Then
$$
\eqalign{
\bigl\|\bigl(e^{-2\pi i\sbmy}-z\id\bigr)v'\bigr\|_\rho^2
&+\bigl\|\bigl(e^{2\pi i\sbmy}-z\id\bigr)\nu'\bigr\|_\rho^2
=\bigl\|\Lambda^{-1}\bigl[e^{2\pi X_0}-z\id\bigr]\Lambda\buv'\bigr\|_\rho^2\cr
&\le\bigl\|\Lambda^{-1}\bigr\|_\rho^2\bigl\|e^{2\pi X_0}\buv-z\buv\bigr\|_\rho^2\,.\cr}
\equation(ResBoundTwo)
$$
Combining this bound with \equ(ResBoundSeven), we get
$$
\eqalign{
\bigl\|\Lambda^{-1}\bigr\|_\rho^2\bigl\|e^{2\pi X_0}\buv-z\buv\bigr\|_\rho^2
&\ge\dist\bigl(z,\Sigma_\rho^e\bigr)^2\|\buv'\|_\rho^2\cr
&\ge\dist\bigl(z,\Sigma_\rho^e\bigr)^2\|\Lambda\|^{-2}\|\buv\|_\rho^2\,.\cr}
\equation(ResBoundEight)
$$
This proves \equ(ResBound).
\qed

Consider the $2n$-dimensional spaces
$\HH_{\rho,n}=\buH_{\rho,n}(\complex)$ defined after \clm(TailBound).
On $\HH_{0,n}$ we define a self-adjoint operator $H_{0,n}$
via the quadratic form
$$
\langle v,H_{0,n}v\rangle
=\bar\alpha|v_{1-n}-v_n|^2+\bar\alpha\sum_{j=1-n}^{n-1}|v_{j+1}-v_j|^2
+\bar\beta\langle v,v\rangle\,,
\equation(Hon)
$$
for $v\in\HH_{0,n}$.
Notice that $H_0=\bar\alpha\nablao^\ast\nablao+\bar\beta$
is invariant under translations.
The operator $H_{0,n}$ is invariant under translations as well,
in the sense that it maps $u$ to $v$ precisely when
$H_0$ maps the $2n$-periodic extension of $u$
to the $2n$-periodic extension of $v$.
Consider also the operator $X_0$ defined in \equ(XoDef).
An operator $X_{0,n}:\HH_{0,n}^2\to\HH_{0,n}^2$ is defined analogously,
with $H_0$ replaced by $H_{0,n}$.

Assume again that $\bar\beta$ and $\bar\beta+4\bar\alpha$ are nonnegative.
Consider the operator $\Lambda$ defined in \equ(XoPartialDiagLambda).
An operator $\Lambda_n:\HH_{0,n}^2\to\HH_{0,n}^2$ is defined analogously,
using $\bmy=H_{0,n}^{1/2}$.

\claim Proposition(ResBoundn)
Let $0\le\rho<\rho_1$.
Then there exists $C>0$ such that for all $n\ge 1$,
$$
C\bigl\|e^{2\pi X_0}\buv-z\buv\bigr\|_\rho
\ge\,\dist\bigl(z,\Sigma_\rho^e\bigr)\|\buv\|_\rho\,,
\qquad\buv\in\HH_{\rho,n}^2\,,\quad z\in\complex\,.
\equation(ResBoundn)
$$

This proposition is proved in the same way as \clm(ResBound).
Uniformity in $n$ follows from translation invariance:
$H_{0,n}$ is diagonalized by the Fourier transform $\FF_n$,
and the Fourier multiplier is always the function $h_0$.

\subsection Spectral approximation

The goal is to approximate the spectrum of $\Flow(2\pi)$
by the spectrum of operators $\Flow_n(2\pi)$ that are essentially matrices.
Define a self-adjoint operator $H_n(t)$ on $\HH_{0,n}$
via the quadratic form
$$
\Langle v,H_n(t)v\Rangle
=\alpha_{1-n}(t)|v_{1-n}-v_n|^2
+\sum_{j=1-n}^{n-1}\alpha_j(t)|v_{j+1}-v_j|^2+\langle v,\beta v\rangle\,,
\equation(Hn)
$$
where $\alpha_j=\tilde\alpha_j+\bar\alpha$
and $\beta_j=\tilde\beta_j+\bar\beta$ for all $j$.
Let us now identify $\HH_{0,n}$ with $P_n\HH_0$,
where $P_n$ is the orthogonal projection defined in \equ(PnDef).
Then $H_{0,n}$ and $H_n$ extend canonically to $\HH_0$
via the identities \equ(Hon) and \equ(Hn), respectively.

In the canonical way we also define
the vector field $X_{0,n}$ associated with $H_{0,n}$,
the time-dependent vector field $X_n$ associated with $H_n$,
and the time-$2\pi$ map $\Flow_n(2\pi)$ associated
with the vector field $X_n$.
Notice that all these operators on $\HH_0^2$ commute with $P_n$
and act trivially on $(\id-P_n)\HH_0^2$.

\claim Proposition(PhinminusPhiV)
Let $0\le\rho<\rho_0$.
Then $\Flow_n(2\pi)\buv\to\Flow(2\pi)\buv$ for every $\buv\in\HH_\rho$.

\proof
Since the operator norms of $\Flow_n(2\pi):\HH_\rho\to\HH_\rho$
are bounded uniformly in $n$,
it suffices to prove that $\Flow_n(2\pi)\buv\to\Flow(2\pi)\buv$
for every $\buv$ in some dense subset of $\HH_\rho$.
To this end, choose $\rho<\varrho<\rho_0$,
and let $\buv$ be a nonzero vector in $\HH_\varrho$.
Define $\buw(t)=\Flow(t)\buv$ and $\buw_n(t)=\Flow_n(t)\buv$.
The difference $\buw-\buw_n$ is the solution of the equation
\vskip-2pt 
$$
{d\over dt}(\buw-\buw_n)=X\buw-X_n\buw_n=X(\buw-\buw_n)+(X-X_n)\buw_n\,,
\equation(WWkOne)
$$
with zero initial condition.
Taking the inner product with $\buw-\buw_n$ yields the bound
$$
\half{d\over dt}\|\buw-\buw_n\|_\rho^2
\le\|X\|_\rho\|\buw-\buw_n\|_\rho^2+\|\buw-\buw_n\|_\rho\|(X-X_n)\buw_n\|_\rho\,,
\equation(WWkTwo)
$$
or equivalently,
\vskip-4pt 
$$
{d\over dt}\|\buw-\buw_n\|_\rho\le\|X\|_\rho\|\buw-\buw_n\|_\rho+\|(X-X_n)\buw_n\|_\rho\,.
\equation(WWkThree)
$$
In order to estimate the last term in this equation,
we use that $\Flow(t)$ and $\Flow_n(t)$
define bounded linear operators on $\HH_\rho^2$,
that $\Flow(t)^{-1}$ and $\Flow_n(t)^{-1}$
define bounded linear operators on $\HH_\varrho^2$,
and that their operator norms are bounded
uniformly in $k$ and in $t\in[0,2\pi]$.
This yields a bound
$$
\|\buw_n\|_\rho\le C_1\|\buv\|_\rho
\le C_1\|\buv\|_\varrho\le C_2\|\buw_n\|_\varrho\,,
\equation(WWkFour)
$$
for some constants $C_1$ and $C_2$.
Here, and in what follows,
a given bound on a quantity that depends on $n$ and $t$
is meant to hold uniformly in $n$ and $t\in[0,2\pi]$.

Using \clm(TailBound), we conclude from \equ(WWkFour) that,
for any given $\eps>0$, there exists $m>0$
such that $\|(\id-P_m)\buw_n\|_\rho\le\eps$ for all $n$.
Furthermore, $\|(X-X_n)P_m\|_\rho\le\eps$ if $n$ is sufficiently large.
Thus, we have a bound
$$
\|(X-X_n)\buw_n\|_\rho
\le\|(X-X_n)(\id-P_m)\buw_n\|_\rho+\|(X-X_n)P_m\buw_n\|_\rho\le\eps_n\,,
\equation(WWkFive)
$$
with $\eps_n\to 0$ as $n\to\infty$.
Now pick $a>\|X\|$. Then \equ(WWkThree) and \equ(WWkFive) imply that
$$
{d\over dt}e^{-at}\|\buw(t)-\buw_n(t)\|_\rho\le\eps_n\,.
\equation(WWkSix)
$$
\vskip-2pt\noindent 
As a result, we have
$$
\|\Flow(2\pi)\buv-\Flow_n(2\pi)\buv\|_\rho
=\|\buw(2\pi)-\buw_n(2\pi)\|_\rho\le\eps_n 2\pi e^{2\pi a}\,.
\equation(WWkSeven)
$$
This proves the claim.
\qed

After having shown that $\Flow_n(2\pi)\to\Flow(2\pi)$ pointwise,
we consider the problem of convergence of eigenvalues and eigenvectors.

\claim Proposition(ApproxEigen)
Let $0<\rho<\rho_0$.
Assume that there exists an increasing sequence $k\mapsto n_k$ of positive integers,
a converging sequence $k\mapsto\lambda_k$ of complex numbers,
and a sequence $k\mapsto\buv_k$ of unit vectors in $\HH_\rho^2$,
such that $\|\Flow_{n_k}(2\pi)\buv_k-\lambda_k\buv_k\|_\rho\to 0$ as $k\to\infty$.
If $\lambda=\lim_k\lambda_k$ does not belong to $\Sigma_\rho^e$,
then $\lambda$ is an eigenvalue of $\Flow(2\pi)$,
and some subsequence of $k\mapsto\buv_k$
converges in $\HH_\rho^2$ to an eigenvector of $\Flow(2\pi)$
for the eigenvalue $\lambda$.

\proof
Assume that $\lambda\not\in\Sigma_\rho^e$.
We may also assume that $\lambda_k\not\in\Sigma_\rho^e$ for all $k$.
Let
$$
E_k=\Flow_{n_k}(2\pi)\buv_k-\lambda_k\buv_k\,.
\equation(AEOne)
$$
Using that $\Flow_{n_k}(2\pi)=e^{2\pi X_{0,n_k}}\bigl[\id+A_k(2\pi)\bigr]$,
we have
$$
\buv_k=\bigl[\lambda_k-e^{2\pi X_{0,n_k}}\bigr]^{-1} e^{2\pi X_{0,n_k}}A_k(2\pi)\buv_k
+\bigl[\lambda_k-e^{2\pi X_{0,n_k}}\bigr]^{-1}E_k\,.
\equation(AETwo)
$$
By \clm(ResBoundn), the operators $\bigl[\lambda_k-e^{2\pi X_{0,n_k}}\bigr]^{-1}$
are bounded on $\HH_\rho^2$, uniformly in $k$.
Thus, the last term in \equ(AETwo) converges to zero as $k\to\infty$.
If $\rho<\varrho<\rho_0$,
then $A_k(2\pi)\buv_k$ belongs to $\HH_\varrho^2$ for each $k$, by \clm(CompactOps).
In fact, it is clear from the proof of \clm(CompactOps)
that the sequence $k\mapsto A_k(2\pi)\buv_k$ is bounded in $\HH_\varrho^2$.
Thus, some subsequence converges in $\HH_\rho^2$.
To simplify notation, we take this subsequence to be $k\mapsto A_k(2\pi)\buv_k$.
Then \equ(AETwo) implies that the sequence $k\mapsto\buv_k$
converges in $\HH_\rho^2$. Denote the limit by $\buv$.
Consider now the inequality
$$
\eqalign{
\|\Flow(2\pi)\buv-\lambda\buv\|_0
&\le\|\Flow(2\pi)\buv-\Flow_{n_k}(2\pi)\buv\|_0
+\|\Flow_{n_k}(2\pi)[\buv-\buv_{n_k}]\|_0\cr
&\quad+\|\Flow_{n_k}(2\pi)\buv_k-\lambda_k\buv_k\|_0
+\|\lambda\buv-\lambda_k\buv_k\|_0\,.\cr}
\equation(Phivminuslambdav)
$$
By \clm(PhinminusPhiV),
the first term on the right hand side tends to zero as $k\to\infty$.
The other three terms on the right tend to zero trivially.
This shows that $\buv$ is an eigenvector of $\Flow(2\pi)$
with eigenvalue $\lambda$.
\qed

\medskip
Denote by $\Sigma_n$ the spectrum of $\Flow_n(2\pi):\HH_0^2\to\HH_0^2$.
It includes one \big(if $\bar\beta=0$\big) or two trivial eigenvalues
$\exp\bigl(\pm i\sqrt{\bar\beta}\,\bigr)$ with infinite multiplicity.
The remaining part of $\Sigma_n$ consists of
eigenvalues whose multiplicities add up to at most $4n$.

\claim Theorem(LimEigen)
Let $\lambda\in\complex\setminus\Sigma_0^e$.
Then $\lambda$ is an eigenvalue of $\Flow(2\pi):\HH_0^2\to\HH_0^2$
if and only if there exists a sequence of points $\lambda_n\in\Sigma_n$
that accumulates at $\lambda$.

\proof
The ``if'' part follows from \clm(ApproxEigen), since we can choose
a positive $\rho<\varrho$ such that $\lambda$ lies outside $\Sigma_\rho^e$.

To prove the ``only if'' part,
assume that $\lambda$ is an eigenvalue of $\Flow(2\pi):\HH_0^2\to\HH_0^2$.
By \clm(EigenDecay), the eigenvectors for this eigenvalue
belong to $\HH_\rho^2$, if $\rho>0$ is chosen sufficiently small.
We may assume that $\rho<\rho_0$.

Choose $r>0$ such that
the closure of the disk $D=\{z\in\complex:|z-\lambda|<r\}$
does not intersect $\Sigma_\rho^e$ and contains no eigenvalue
of $\Flow(2\pi)$ besides $\lambda$.
Then the spectral projection associated
with $D$ of the operator $A=\Flow(2\pi)$ is given by
$$
\proj(A,D)={1\over 2\pi i}\int_{\partial D}(z\id-A)^{-1}\,dz\,,
\equation(RieszProj)
$$
where $\partial D$ denotes the positively oriented boundary of $D$.
The goal is to show that the corresponding projection
for $A=\Flow_n(2\pi)$ is well defined and nontrivial,
if $n$ is sufficiently large.
To this end, define
$$
c=\liminf_{n\to\infty}\,\inf_{z\in\partial D}\,\inf_{\sbuv\in B}\|\Flow_n(2\pi)\buv-z\buv\|_\rho\,,
\equation(LEOne)
$$
where $B$ is the unit ball in $\HH_\rho^2$.
Assume for contradiction that $c=0$.
Then we can find an increasing sequence $k\mapsto n_k$,
a sequence $k\mapsto\lambda_k\in\partial D$,
and a sequence $k\mapsto\buv_k\in B$,
such that $\|\Flow_{n_k}(2\pi)\buv-\lambda_k\buv\|_\rho\to 0$.
By choosing a subsequence, if necessary,
we can achieve $\lambda_k\to\lambda\in\partial D$.
By \clm(ApproxEigen), this implies that $\lambda$
is an eigenvalue of $\Flow(2\pi)$.
But this is impossible by our choice of $D$.
The conclusion is that $c>0$.

Thus, there exists $N>0$ that
$$
\|\Flow_n(2\pi)\buv-z\buv\|_\rho\ge{c\over 2}\|\buv\|_\rho\,,
\qquad z\in\partial D\,,\quad n\ge N\,,
\equation(LETwo)
$$
for all $\buv\in\HH_\rho^2$.
This implies e.g.~that the spectral projection \equ(RieszProj)
is well-defined for $A=\Flow_n(2\pi)$.
Here, and in what follows, we assume that $n\ge N$.
Let now $\buv\in\HH_\rho$ be an eigenvector of $\Flow(2\pi)$
with eigenvalue $\lambda$.
By the second resolvent identity we have
$$
\Langle\buv,\bigl[
\proj(\Flow_n(2\pi),D)-\proj(\Flow(2\pi),D)\bigr]\buv\Rangle_0
={1\over 2\pi i}\int_{\partial D}f_n(z)\,dz\,,
\equation(LEThree)
$$
where
$$
f_n(z)=\Langle\buv,\bigl(z\id-\Flow_n(2\pi)\bigr)^{-1}E_n(z)\Rangle_0
\equation(LEFour)
$$
and
$$
E_n(z)=\bigl[\Flow_n(2\pi)-\Flow(2\pi)\bigr]\bigl(z\id-\Flow(2\pi)\bigr)^{-1}\buv\,.
\equation(LEFive)
$$
The goal is to take $n\to\infty$.
The bound \equ(LETwo) implies that the sequence $n\mapsto f_n$ is bounded
uniformly on $\partial D$.
Furthermore, we have $\|E_n(z)\|_0\to 0$ for each $z\in\partial D$,
as a consequence of \clm(PhinminusPhiV).
By \clm(ResBoundn), this implies that $f_n\to 0$ pointwise on $\partial D$.
And by the bounded convergence theorem,
it follows that the integral in \equ(LEThree) tends to zero as $n\to\infty$.

Given that
$\langle\buv,\proj\bigl(\Flow(2\pi),D\bigr)\buv\rangle_0=\|\buv\|_0^2>0$,
we conclude that
$\langle\buv,\proj\bigl(\Flow_n(2\pi),D\bigr)\buv\rangle_0$
is nonzero for sufficiently large $n$.
This shows that $\Flow_n(2\pi)$ has an eigenvalue in $D$,
if $n$ is sufficiently large.
Since the radius $r>0$ of $D$ can be taken arbitrarily small,
the assertion follows.
\qed

\section Existence of solutions

\vskip-15pt

\subsection Localized solutions

In this section we prove \clm(breathers),
based on a technical lemma that will be proved later.

Adding $\psi_2 q$ on both sides of the equation \equ(Main),
we obtain
$$
\bigl(\omega^2\partial_t^2+\psi_2\bigr)q
=-\nablas^\ast\phi'(\nablas q)-{\tilde\psi}'(q)\,,\qquad
\tilde\psi(x)=\psi(x)-{\textstyle{1\over 2}}\psi_2 x^2\,.
\equation(preGDef)
$$
Formally, we can rewrite \equ(preGDef)
as the fixed point equation $G(q)=q$, where
$$
G(q)=-L^{-1}\bigl[\nablas\phi'(\nablas q)+{\tilde\psi}'(q)\bigr]\,,\qquad
L=\omega^2\partial_t^2+\psi_2\,.
\equation(GDef)
$$
For the domain of $G$ we use one of the spaces
$\BB_{r,\rho}^\ssigma$ defined below.

Given a real number $r>1$, denote by $\AA_r$ the Banach space
of all $2\pi$-periodic functions $g:\real\to\real$
that have a finite norm $\|g\|_r$,
$$
g(t)=\sum_{k\ge 0}g_k^\even\cos(kt)
+\sum_{k\ge 1}g_k^\odd\sin(kt)\,,\qquad
\|g\|_r=\sum_{k\ge 0}\bigl|g_k^\even\bigr|r^k
+\sum_{k\ge 1}\bigl|g_k^\odd\bigr|r^k\,.
\equation(CosSinSeries)
$$
The even and odd subspaces of $\AA_r$
are denoted by $\AA_r^\even$ and $\AA_r^\odd$, respectively.
A straightforward computation shows that
$\|fg\|_r\le\|f\|_r\|g\|_r$.
Thus, $\AA_r$ and $\AA_r^\even$ are Banach algebras
under pointwise multiplication.
Notice that the functions in $\AA_r$ extend analytically
to the complex domain $|\Im(t)|<\log(r)$.

Next, consider the vector space of all chains
$j\mapsto v_j\in\AA_r$ with only finitely many nonzero values $v_j$.
Such chains will be called finite.
Given $\sigma\in\{0,1\}$ and a real number $\rho>1$,
we define $\BB_{r,\rho}^\ssigma$ to be the completion of this space
with respect to the norm
$$
\|v\|_{r,\rho}^\sigma=\sum_{j\in\integer}\|v_j\|_r\rho^{|2j-\sigma|}\,.
\equation(ChainNorm)
$$
The even and odd (as functions of $t$) subspaces of $\BB_{r,\rho}^\ssigma$
are denoted by $\BB_{r,\rho}^{\ssigma,\even}$ and $\BB_{r,\rho}^{\ssigma,\odd}$, respectively.
A straightforward computation shows that
$$
\|uv\|_{r,\rho}^\sigma\le\|u\|_{r,\rho}^\sigma\sup_j\|v_j\|_r\,,\qquad
\|v_j\|_r\le\rho^{-|2j-\sigma|}\|v\|_{r,\rho}^\sigma\,.
\equation(uvNorm)
$$
In particular, $\BB_{r,\rho}^\ssigma$ and $\BB_{r,\rho}^{\ssigma,\even}$
are Banach algebras under pointwise multiplication $(uv)_j=u_j v_j$.
Furthermore, chains in $\BB_{r,\rho}^\ssigma$ decrease exponentially.
The operator norm of a continuous linear operator $L$
on $\BB_{r,\rho}^\ssigma$ or $\BBB$ will be denoted by $\|L\|_{r,\rho}^\sigma$.

Notice that $\BB_{r,\rho}^{\,0}=\BB_{r,\rho}^{\,1}$ as vector spaces.
Our reason for choosing the $\sigma$-dependent norm \equ(ChainNorm)
is that the reflection $q\mapsto\breve q$ defined by \equ(symmetry)
is an isometry for this norm.

\smallskip
In order to solve the fixed point problem for $G$,
we first determine (numerically) a finite chain $\bar q$
that is an approximate fixed point of $G$,
and a linear isomorphism $A$ of $\BB_{r,\rho}$
that is an approximate inverse of $\id-DG(\bar q)$.
Then the map $\FF$ defined by
$$
\FF(h)=G(q)-q+h\,,\qquad q=\bar q+Ah\,,
\equation(contr)
$$
can be expected to be a contraction near the origin.
Clearly, $h$ is a fixed point of $\FF$
if and only if $q$ is a fixed point of $G$.

Consider now a fixed but arbitrary row in Table 1.
Among other things, it specifies
a domain parameter $\sBB=(\sigma,\pars,r,\rho)$ identifying a space $\BBB$,
and a symmetry parameter $\varpi$.
If $\varpi=\pm$, then we define $\RR_{\sigma,\varpi}$ to be the reflection
$q\mapsto\breve q$ given by the equation \equ(symmetry).
Otherwise, if $\varpi=$``none'',
then $\RR_{\sigma,\varpi}$ is defined to be the identity map.

The following lemma is proved with the assistance
of a computer, as described in Section 5.

\claim Lemma(contraction)
For each set of parameters $(\phi,\psi,\omega,\varpi,\sBB,\ell)$ given in Table 1,
there exists a finite chain $\bar q$,
a linear isomorphism $A$ of $\BBB$,
and positive constants $\eps,K,\delta$ satisfying $\eps+K\delta<\delta$,
such that the map $\FF$ defined by \equ(contr)
is analytic on $\BBB$ and satisfies
$$
\|\FF(0)\|_{r,\rho}^\sigma\le\eps\,,\qquad
\|D\FF(h)\|_{r,\rho}^\sigma\le K\,,\qquad h\in B_{\delta,2\delta}\,,
\equation(contraction)
$$
with $B_{\delta,2\delta}$ as defined below.
The support of $\bar q$ is the set $\{j\in\integer:\sigma-\ell<j<\ell\}$,
and $Ah=h$ for every chain $h$ that vanishes on this set.
Furthermore, $\bar q$ is invariant under $\RR_{\sigma,\varpi}$
and $A$ commutes with $\RR_{\sigma,\varpi}$.

\smallskip
The set $B_{\delta,2\delta}$ in the equation \equ(contraction)
is a special case of the following.
Given real numbers $u,v>0$, we define
$B_{u,v}$ to be the set of all chains $h\in\BBB$ with the property that
$$
\sum_{j\le\sigma-\ell}\|h_j\|_r\rho^{|2j-\sigma|}\le v\,,\quad
\!\!\sum_{\sigma-\ell<j<\ell}\!\!\|h_j\|_r\rho^{|2j-\sigma|}\le u\,,\quad
\sum_{j\ge\ell}\|h_j\|_r\rho^{|2j-\sigma|}\le v\,.
\equation(Bdelta)
$$
Notice that $B_{\delta,2\delta}$ includes the closed ball in $\BBB$
of radius $\delta$, centered at the origin.

\proofof(breathers)
First, we note that $\phi$ and $\psi$ are even whenever $\pars=1$.
Thus, the right hand side of \equ(preGDef)
belongs to $\BB=\BBB$ whenever $q\in\BB$.
Furthermore, $\psi_2-(\omega k)^2$ is bounded away from zero for all integers $k\ge\pars$.
This implies that $L:\BB\to\BB$ has a bounded inverse.
Thus, $G$ is well-defined on $\BB$ and analytic (in fact polynomial).
The same is true for the map $\FF$, since $A$ is bounded.

By the contraction mapping theorem, the given bounds imply that
$\FF$ has a unique fixed point $h_\ast\in\BBB$ with norm $\le\delta$.
Now $q_\ast=\bar q+Ah_\ast$ is a fixed point of $G$
and thus satisfies the equation \equ(Main).

It is straightforward to check that $G$ commutes with $\RR=\RR_{\sigma,\varpi}$.
Since $A$ commutes with $\RR$ as well,
the same is true for the map $\FF$.
Here, we have used also that $\bar q$ is invariant under $\RR$.
Thus, given that $\displaystyle h_\ast=\lim_{k\to\infty}\FF^k(0)$,
it follows that $h_\ast$ and $q_\ast=\bar q+Ah_\ast$
are invariant under $\RR$.
\qed

\medskip
We note that an alternative to the map $G$ considered here
would be the map $\GG$, defined by
$$
\GG(q)=-\LL^{-1}\bigl[\nablas{\tilde\phi}'(\nablas q)+{\tilde\psi}''(q)\bigr]\,,\qquad
\LL=\omega^2\partial_t^2+\phi_2\nablas^\ast\nablas+\psi_2\,,
\equation(GGDef)
$$
where $\tilde\phi(x)=\phi(x)-\half\phi_2x^2$.
The inverse of $\LL$ involves a lattice-convolution
with an exponentially decreasing kernel
(for suitable values of $\phi_2$, $\psi_2$, and $\omega$).
This kernel can be computed explicitly,
but its nonlocality complicates the analysis significantly;
especially the construction of multi-breather solutions.

\subsection Combining solutions

In this subsection we give a general result
that will be used later to prove \clm(multi).
This part is independent of previous sections,
which allows us to adapt the notation to the problem at hand.

In what follows, if we write a Banach space $Y$
as a direct sum of subspaces,
$$
Y=\bigoplus_k Y_k\,,
\equation(YYkDirSum)
$$
then the norm on this space is assumed to satisfy
$$
\|y\|_{\sss Y}=\sup_k\|y_k\|_{\sss Y}\,,\qquad y_k=P_k y\,,
\equation(YYkDirSumNorm)
$$
where $P_k$ denotes the canonical projection from $Y$ onto $Y_k$.
Let now $Y$ be a direct sum as in \equ(YYkDirSum),
where the index $k$ runs over the set of all integers.

For each integer $k$,
let $f_k$ be a $\rmC^1$ mapping on $Y_k\oplus Y_{k+1}$.
We extend $f_k$ to $Y$ by setting $f_k(y)=f_k(y_k,y_{k+1})$.
Define
$$
F_n^{-}(y)={\sum_{k<n}}'f_k(y)\,,\qquad
F_n^{+}(y)={\sum_{k\ge n}}'f_k(y)\,,\qquad
F(y)={\sum_k}'f_k(y)\,,
\equation(FnDef)
$$
for all $y\in Y$.

\claim Notation(Pointwise)
Here, and in what follows, $\sum'$ denotes a pointwise sum,
meaning that its $m$-th component converges in $Y_m$, for each $m$.

Notice that, for the sums in \equ(FnDef),
each component $P_m\sum_k'f_k(y)$ is a sum of at most two nonzero terms.

For odd integers $n$,
let $X_n\subset\bigoplus_{k\le n}Y_k$ and $Z_n\subset\bigoplus_{k\ge n}Y_k$
be subspaces that carry norms $\|\bdot\|_{\sss X_n}$
and $\|\bdot\|_{\sss Z_n}$, respectively,
and that are complete for these norms.
We also assume that $X_n\cap Z_n=Y_n$ and
$$
\|y_n\|_{\sss X_n}=\|y_n\|_{\sss Z_n}=\|y_n\|_{\sss Y}\,,\qquad y_n\in Y_n\qquad
(n {\rm ~odd}).
\equation(XZNormCond)
$$

In the remaining part of this subsection, we always assume that
$n$ is an \ub{even} integer, unless specified otherwise.

Define $\HH_n=X_{n-1}\oplus Y_n\oplus Z_{n+1}$.
So a vector $h\in\HH_n$ admits a unique representation $h=x_{n-1}+y_n+z_{n+1}$
with $x_{n-1}\in X_{n-1}$, $y_n\in Y_n$, and $z_{n+1}\in Z_{n+1}$.
We will also use the notation $h=(x_{n-1},y_n,z_{n+1})$.

\demo Remark(XYZmeaning)
We think of $h=(x_{n-1},y_n,z_{n+1})$ as representing a chain
with ``center'' $y_n$, left-tail $x_{n-1}$, and right-tail $z_{n+1}$.
The idea is to take $F$ locally of the form \equ(contr),
with $\bar q$ and $A$ depending on $n$.
The goal is to show that $F$ has a fixed point near the origin in $Y$.

We now impose conditions on the function $F$
that can be checked separately for each of the spaces $\HH_n$.
Let $\eps,K,\delta$ be positive integers
satisfying $\eps+K\delta<\delta$.
Assume that
$$
\|P_nF(0)\|_{\sss Y_n}\le\eps\,,\quad
\|P_{n-1}F_n^{-}(0)\|_{\sss X_{n-1}}\le\thalf\eps\,,\quad
\|P_{n+1}F_n^{+}(0)\|_{\sss Z_{n+1}}\le\thalf\eps\,.
\equation(ApproxSolution)
$$
In addition, assume
that $F_n^{-}$ defines a $\rmC^1$ function on $X_{n-1}\oplus Y_n$,
that $F_n^{+}$ defines a $\rmC^1$ function on $Y_n\oplus Z_{n+1}$,
and that
$$
\|P_nDF(h)\|_{\sss\HH_n\to Y_n}\le K\,,
\equation(LocalContraction)
$$
whenever $\|h\|_{\sss\HH_n}\le\delta$.

In order to formulate our last assumption,
we write $F_n^{-}$ as a function of two arguments,
the first in $X_{n-1}$ and the second in $Y_n$.
Similarly, we write $F_n^{+}$ as a function of two arguments,
the first in $Y_n$, and the second in $Z_{n+1}$.
Assume that
$$
\eqalign{
\|P_{n-1}DF_n^{-}(x_{n-1},y_n)\|_{\sss X_{n-1}\oplus Y_n\to Y_{n-1}}&\le\thalf K\,,\cr
\|P_{n+1}DF_n^{+}(y_n,z_{n+1})\|_{\sss Y_n\oplus Z_{n+1}\to Y_{n+1}}&\le\thalf K\,,\cr}
\equation(ErrDer)
$$
whenever $\|h\|_{\sss\HH_n}\le\delta$.

\claim Proposition(yCombined)
Under the assumptions described above,
which includes the condition $\eps+K\delta<\delta$,
the map $F:Y\to Y$ has a unique fixed point
in the closed ball of radius $\delta$ in $Y$,
centered at the origin.

\proof
If $n$ is odd, then
$$
\|P_nF(0)\|_{\sss Y_n}\le\|P_n F_{n-1}^{+}(0)\|_{\sss Y_n}
+\|P_n F_{n+1}^{-}(0)\|_{\sss Y_n}\le\thalf\eps+\thalf\eps\,,
\equation(ycOne)
$$
by the second and third inequality in \equ(ApproxSolution).
Combining this with the first inequality in \equ(ApproxSolution) yields
$$
\|F(0)\|_{\sss Y}=\sup_{n\in\integer}\bigl\|P_nF(0)\|_{\sss Y_n}\le\eps\,.
\equation(yCTwo)
$$
Let now $u$ be a fixed but arbitrary vector in $Y$ with norm $\|u\|\le\delta$.
Our goal is to estimate $DF(u)$.
Notice that
$$
\eqalign{
P_nDF(u)v
&=P_nD_{\sss Y_{n-1}}f_{n-1}(u_{n-1},u_n)v_{n-1}
+P_nD_{\sss Y_n}f_{n-1}(u_{n-1},u_n)v_n\cr
&\quad+P_nD_{\sss Y_n}f_n(u_n,u_{n+1})v_n
+P_nD_{\sss Y_{n+1}}f_n(u_n,u_{n+1})v_{n+1}\,,\cr}
\equation(yCThree)
$$
for all $u,v\in Y$ and all integers $n$.
Here, $D_{\sss Y_k}$ denote the partial derivative operator with respect
to the component in $Y_k$.

Consider first the case where $n$ is even.
Setting $x_{n-1}=u_{n-1}$, $y_n=u_n$, and $z_{n+1}=u_{n+1}$,
the vector $h=(x_{n-1},y_n,z_{n+1})$ has norm $\|h\|_{\HH_n}\le\delta$.
So by \equ(yCThree) and \equ(LocalContraction), we have
$$
\|P_nDF(u)\|_{\sss Y\to Y_n}
\le\|P_nDF(x_{n-1},y_n,z_{n+1})\|_{\sss\HH_n\to Y_n}\le K\,.
\equation(yCFour)
$$
Next, consider the case when $n$ is odd.
Setting  $y_{n-1}=u_{n-1}$, $x_n=z_n=u_n$, and $y_{n+1}=u_{n+1}$, we have
$$
\eqalign{
\|P_nDF(u)\|_{\sss Y\to Y_n}
&\le\|P_nDF_{n-1}^{+}(y_{n-1},z_n)\|_{\sss Y_{n-1}\oplus Z_n\to Y_n}\cr
&\quad+\|P_nDF_{n+1}^{-}(x_n,y_{n+1})\|_{\sss X_n\oplus Y_{n+1}\to Y_n}
\le\thalf K+\thalf K\,,\cr}
\equation(yCFive)
$$
by \equ(yCThree) and \equ(ErrDer).
Combining \equ(yCFour) and \equ(yCFive) yields
$$
\|DF(u)\|_{\sss Y\to Y}
=\sup_{n\in\integer}\|P_nDF(u)\|_{\sss Y\to Y_n}\le K\,.
\equation(yCSix)
$$
The claim now follows from the contraction mapping theorem.
\qed

\subsection Multi-breather solutions

In this subsection, we consider a fixed but arbitrary
choice of parameters $(\phi,\psi,\omega,\pars,\sigma,r,\rho)$
that is represented by one of the rows $1$-$4$ or $6$-$12$ of Table 1.
The claim in \clm(multi) is that we can produce solutions
that look like strings of breather solutions.
To simplify notation, consider first the case of a bi-infinite string,
indexed by $\integer$.

For any given integer $m$, we choose one of the maps $\FF$
for the given parameters, as described in \clm(contraction).
This involves an approximate fixed point $\bar q$ of $G$ and an operator $A$.
We note that, if $q$ is a possible choice
for the finite chain $\bar q$ mentioned in \clm(contraction),
then $-q$ is an equally good choice.
Here we allow either choice.

Let $A'=A-\id$. After choosing an integer $J_m$,
we set $\bar q^{\tinyskip m}=\TT_{J_m}\bar q$ and $A_m'=T_{J_m}A'T_{J_m}^{-1}$.
Here $T_J$ denotes translation by $J$,
that is, $(T_Jh)_j=h_{j-J}$ for all $j$.
Using the positive integer $\ell$ from Table 1,
define $\JJ_m=\{j\in\integer:\sigma-\ell\le j-J_m\le\ell\}$.
This is the set that we referred to as the approximate support
of the breather $q^m$ in \clm(multi).
It includes the support of $\bar q^{\tinyskip m}$.

We may assume that the sequence $m\mapsto J_m$ is increasing, and that $J_0=0$.
Assuming furthermore that $\dist(\JJ_m,\JJ_{m+1})$
is positive and even for all $m$, we define
$$
F(h)=G(h)+{\sum_m}'\Bigl[
G(h+\bar q^{\tinyskip m}+A_m'h)-G(h)-\bar q^{\tinyskip m}-A_m'h\Bigr]\,.
\equation(multiF)
$$
The sum in this equation converges pointwise, at each integer $j\in\integer$,
since $(A_m'h)_j=0$ whenever $j$ lies outside the support of $\bar q^{\tinyskip m}$.

In order to see how this fits into the framework
discussed in the preceding subsection,
consider a fixed term in this sum, indexed by $m$.
Let $n=2m$.
Consider the translated space $\BB_n=T_{J_m}\BBB$
with norm $\|h\|_{\BB_n}=\bigl\|T_{J_m}^{-1}h\bigr\|_{r,\rho}^\sigma$.
To every chain $h\in\BB_n$
we associate its left-tail $x\in\BB_n$, center $y\in\BB_n$,
and right-tail $z\in\BB_n$ by setting
$$
x_j=\cases{h_j &if $j\le j_m^{-}$,\cr 0 &if $j>j_m^{-}$,\cr}\quad
y_j=\cases{h_j &if $j_m^{-}<j<j_m^{+}$,\cr 0 &otherwise,\cr}\quad
z_j=\cases{h_j &if $j\ge j_m^{+}$,\cr 0 &if $j<j_m^{+}$,\cr}
\equation(xjyjzj)
$$
where $j_m^{-}=\half[\max(\JJ_{m-1})+\min(\JJ_m)]$
and $j_m^{+}=\half[\max(\JJ_m)+\min(\JJ_{m+1})]$.
Denote by $P_n$ the projection $h\mapsto y$ and set $Y_n=P_n\BB_n$.
The ranges of the projections $h\mapsto x$ and $h\mapsto z$
are denoted by $X_{n-1}$ and $Z_{n+1}$, respectively.
In addition, we define $Y_{n-1}$ and $Y_{n+1}$ to be the one-dimensional subspaces of $\BB_n$
spanned by all chains supported in $\{j_m^{-}\}$ and $\{j_m^{+}\}$, respectively.
On $\HH_n=X_{n-1}\oplus Y_n\oplus Z_{n-1}$ we choose the norm
$$
\|h\|_{\HH_n}=
\max\bigl\{\|x\|_{\sss\BB_n}\,\bcomma\,\|y\|_{\sss\BB_n}\,\bcomma\,\|z\|_{\sss\BB_n}\bigr\}\,.
\equation(HnNorm)
$$
The goal now is to apply \clm(yCombined).
The following is meant to be a continuation of \clm(contraction).

\claim Lemma(bumpcontr)
Consider one of the rows $1$-$4$ or $6$-$12$ in Table 1.
In addition to the properties of $\FF$ described in \clm(contraction),
the bounds \equ(LocalContraction) and \equ(ErrDer) with $n=0$
are satisfied for each $h\in B_{\delta,2\delta}$.
Furthermore, $\delta<2^{-50}$ and  $\delta\|A\|_{r,\rho}^\sigma<2^{-46}$.

For the proof of these estimates, we refer to Section 5.
Based on this lemma, we can now give a

\proofof(multi).
Consider first the case of a bi-infinite string.
Then we may assume that the index set is $\integer$,
and that the sequence $m\mapsto J_m$ has the properties
mentioned before \equ(multiF).
The goal is to verify the assumptions of \clm(yCombined).
The parameters $(\phi,\psi,\omega,\sigma,\pars,r,\rho)$ are assumed to be fixed.

By translation invariance, it suffices to verify the conditions
\equ(ApproxSolution), \equ(LocalContraction), and \equ(ErrDer) for $n=0$.
Due to the projections that appear in these conditions,
$F$ can be replaced by the map $\FF$ associated with $\bar q=\bar q^{\tinyskip 0}$ and $A=A_0$.
Notice that the set $B_{\delta,2\delta}$ defined by \equ(Bdelta)
includes the closed ball in $\HH_0$ of radius $\delta$, centered at the origin.
Thus, under our assumption that \equ(LocalContraction) and \equ(ErrDer)
hold for $h\in B_{\delta,2\delta}$,
these bounds hold whenever $\|h\|_{\HH_0}\le\delta$,
as required by \clm(yCombined).
The first inequality in \equ(ApproxSolution)
follow from the first inequality in \equ(contraction).
The other two inequalities in \equ(ApproxSolution)
are satisfied trivially in our case: $j_0^{\pm}$ is at a distance $\ge 2$
from the support of any of the chains $\bar q^{\tinyskip m}$,
so $P_{\pm 1}F_0^{\pm}(0)=0$.

\clm(yCombined) now implies that $F$ has a locally unique fixed point $h\in Y$.
Clearly, the chain
$$
q=h+{\sum_m}'\bigl[\bar q^{\tinyskip m}+A_m' h\bigr]
\equation(qmulti)
$$
is a solution of the equation \equ(Main). Here we are using \clm(Pointwise).

Notice that $\bar q^{\tinyskip m}+A_m' h$ is supported in $\JJ_m$, for each $m$.
For $j$ in between those supports,
we have $|q_j(t)|=|h_j(t)|\le\|h\|_Y\le\delta$.
Consider now  $j\in\JJ_0$.
If $h^0$ denoted the fixed point of the map
$\FF$ associated with $\bar q=\bar q^{\tinyskip 0}$ and $A=A_0$,
and if $q^0=\bar q+Ah^0$ denotes the corresponding solution of \equ(Mainj),
then
$$
\bigl|q_j-q^0_j\bigr|=\bigl|\bigl(A\bigl(h-h^0\bigr)\bigr)_j\bigr|
\le\|A\|_{\sss Y_0\to Y_0}\bigl\|P_0\bigl(h-h^0\bigr)\bigr\|_{\sss Y_0}
\le 2\delta\|A\|_{r,\rho}^\sigma\,,
\equation(qminusqastj)
$$
for all $j\in\JJ_0$.
The same bound holds of course for $j\in\JJ_m$ and any $m$.
This concludes the proof of \clm(multi)
for the case of two-sided infinite strings $m\mapsto\bar q^{\tinyskip m}$.
The proof for one-sided infinite strings
and for finite strings is similar, so we omit it here.
\qed

\section Spectral estimates

Our goal is to reduce the proof of \clm(stability)
to estimates on finite-dimensional systems.

\subsection Instability

We first consider the task of proving spectral instability.
Let $\HH=\ell^2(\integer)$.
The simplest cases are the solutions $5$, $6$, $7$, and $12$,
where the set $\Sigma_0^e$ defined by \equ(EssSpec)
includes a real interval containing the point $1$.
These solutions are spectrally unstable,
since $\Sigma_0^e$ is the essential spectrum of $\Flow(2\pi):\HH^2\to\HH^2$,
as was described before \equ(EssSpec).

In the other cases we use a perturbation argument,
involving an approximation $\Flow_o(2\pi)$ for the map $\Flow(2\pi)$.
First, we need a uniform bound.
Let $(j,t)\mapsto\alpha_j(t)$ and $(j,t)\mapsto\beta_j(t)$
be bounded functions on $\integer\times\real$
that are continuous in the time variable $t$.
Consider the flow on $\HH^2$ given by the equation
$$
\partial_t\buv=X\buv\,,\qquad\buv=\stwovec{v}{\nu}\,,
\quad X=\stwomat{0}{\id}{-H}{0}\,,\quad
Hv=\nabla_\sigma^\ast\alpha\nabla_\sigma v+\beta v\,.
\equation(vFlowAgain)
$$
Here $\alpha v$ and $\beta v$ are defined by pointwise multiplication.
Assume that we have enclosures
$[\alpha_\ast^{-},\alpha_\ast^{+}]\ni\alpha_j(t)$
and $[\beta_\ast^{-},\beta_\ast^{+}]\ni\beta_j(t)$
that are valid for all $j$ and all $t$.
Define
$$
c=\thalf\max\bigl\{|4\alpha_\ast^{-}+\beta_\ast^{-}-1|
\,\bcomma\,|4\alpha_\ast^{+}+\beta_\ast^{+}-1|\bigr\}\,.
\equation(cDef)
$$

\claim Proposition(PhitBound)
Under the above-mentioned assumption,
the time-$t$ maps $\Flow(t)$ for the flow \equ(vFlowAgain)
satisfy the bounds $\|\Flow(t)\|\le e^{c|t|}$ for all $t$.

\proof
Let $\buv=\buv(t)$ be a fixed but arbitrary solution of the equation \equ(vFlowAgain).
Using that $\|H-\id\|\le 2c$, we have
$$
\eqalign{
\partial_t\|\buv\|^2
&=\Langle\buv,(X+X^\ast)\buv\Rangle
=\Langle v,(\id-H)\nu\Rangle+\Langle\nu,(\id-H)v\Rangle\cr
&\le\|\id-H\|\|\buv\|^2\le 2c\|\buv\|^2\,,\cr}
\equation(PhitBoundOne))
$$
and thus $\partial_t\|\buv\|\le c\|\buv\|$,
for every $t\in\real$. By integration we obtain
$$
\bigl\|\Flow(t)\buv(0)\bigr\|
=\|\buv(t)\|\le e^{c|t|}\|\buv(0)\|\,.
\equation(PhitBoundTwo)
$$
This holds for arbitrary initial conditions $\buv(0)\in\HH^2$.
Thus, $\|\Flow(t)\|\le e^{c|t|}$ as claimed.
\qed

For simplicity, assume now that $\tilde\alpha=\alpha-\bar\alpha$
and $\tilde\beta=\beta-\bar\beta$ satisfy the \clm(alphabetaProp),
with $\bar\beta$ and $\bar\beta+4\bar\alpha$ contained in $[0,1)$.
Then the spectrum of $\Flow(2\pi)$ off the unit circle
consists of isolated eigenvalues with finite multiplicities.

Consider another operator $H^o$ of the same type,
for sequences $\alpha^o$ and $\beta^o$
that have the same asymptotic limits $\bar\alpha$ and $\bar\beta$.

\claim Proposition(unstableBound)
Let $\alpha_\ast^{\pm}$ and $\beta_\ast^{\pm}$ be real numbers,
such that $\{\alpha^o_j(t),\alpha_j(t)\}\subset [\alpha_\ast^{-},\alpha_\ast^{+}]$
and $\{\beta^o_j(t),\beta_j(t)\}\subset [\beta_\ast^{-},\beta_\ast^{+}]$
holds for all $j$ and all $t$.
Let $\mu$ be a real number larger than $1$.
Assume that the time-$2\pi$ map $\Flow_o(2\pi)$
associated with $H^o$ has an odd number of eigenvalues
(counting multiplicities) in the half-plane $\Re(z)>\mu$,
and that $\mu$ is not an eigenvalue of $\Flow_o(2\pi)$.
If in addition,
$$
2\pi e^{2\pi c}\bigl\|H(t)-H^o(t)\bigr\|
\bigl\|(\Flow_o(2\pi)-\mu)^{-1}\bigr\|<1
\equation(uBoundOne)
$$
for all $t\in[0,2\pi]$, with $c$ given by \equ(cDef),
then $\Flow(2\pi)$ has an odd number of eigenvalues
in the half-plane $\Re(z)>\mu$.

\proof
For $0\le\kappa\le 1$ define $H^\kappa=(1-\kappa)H^o+\kappa H$.
Denote by $A_\kappa$ the time-$2\pi$ map associated with $H^\kappa$.
Our goal is to show that
$$
\bigl\|(A_\kappa-A_0)(A_0-\mu)^{-1}\bigr\|<1\,,\qquad
0\le\kappa\le 1\,.
\equation(uBoundTwo)
$$
Then each $A_\kappa-\mu=\bigl[\id+(A_\kappa-A_0)(A_0-\mu)^{-1}\bigr](A_0-\mu)$
has a bounded inverse, implying that no $A_\kappa$
has an eigenvalue $\mu$.
Since the eigenvalues of $A_\kappa$ off the unit circle
depend continuously on $\kappa$ and come in complex-conjugate pairs,
this implies that each operator $A_\kappa$ has an odd number of eigenvalues
in the half-plane $\Re(z)>\mu$.
So the claim made in \clm(unstableBound) follows form the bound \equ(uBoundTwo).

What we need now is a bound on $A_\kappa-A_0$.
Denote by $\Flow^\kappa_{t,s}$ the flow-map for $H^\kappa$ from time $s$ to time $t$.
These maps satisfy the equation
$$
\Flow^\kappa_{t,r}=\Flow^0_{t,r}+\int_r^t
\Flow^0_{t,s}P_2^\ast\bigl[H^\kappa(s)-H^o(s)\bigr]P_1\Flow^\kappa_{s,r}\,ds\,,
\equation(uBoundThree)
$$
where $P_1=[1\ 0]$ and $P_2=[0\ 1]$ are the
the operators from $\HH^2$ to $\HH$ that are described after \equ(PhiIntEqu).
By \clm(PhitBound) we have $\|\Flow^\kappa_{t,s}\|\le e^{c|t-s|}$ for each $\kappa$.
Taking norms in \equ(uBoundThree) we get
$$
\bigl\|\Flow^\kappa_{t,r}-\Flow^0_{t,r}\bigr\|
\le\int_r^t e^{c(t-s)}\bigl\|H^\kappa(s)-H^o(s)\bigr\|e^{c(s-r)}\,ds\,.
\equation(Dtwo)
$$
In particular,
$$
\|A_\kappa-A_0\|\le 2\pi\kappa e^{2\pi c}\sup_{0\le t\le 2\pi}
\bigl\|H(t)-H^o(t)\bigr\|\,.
\equation(Dthree)
$$
When combined with the assumption \equ(uBoundOne),
this yields the desired bound \equ(uBoundTwo).
\qed

\smallskip
Our choice of $H^o$ will be described in Subsection 4.4.

\subsection Separating sets and monotonicity

Next, we consider the task of proving spectral stability.
We adapt an approach that was introduced in [\rAKii].
Roughly speaking, the goal is to find two simple
approximations $H^{\pm 1}$ for the operator $H$ defined in \equ(HqDef),
such that $H^{-1}\ltstrong H\ltstrong H^1$.
If we can control the time-$2\pi$ maps $\Flow_s(2\pi)$
associated with the family of operators $H^s={1-s\over 2}H^{-1}+{1+s\over 2}H^1$,
in a way that will be explained below,
then we can also control the time-$2\pi$ map $\Flow(2\pi)$ associated with $H$.
And as described at the end of this subsection,
we can reduce this to a finite-dimensional problem.
Here, and in what follows, the parameter $s$ always ranges
over the interval $[-1,1]$.

Let $\HH$ be a finite-dimensional Hilbert space.
Consider the Hilbert space $\HH^2$ of all pairs $\buv=\bigl[{v\atop\nu}\bigr]$,
equipped with the inner product
$\langle\buv,\buv'\rangle=\langle v,v'\rangle+\langle\nu,\nu'\rangle$.
A linear operator on $\HH^2$ is said to be symplectic
if it is ``unitary'' for the quadratic form
$$
G(\buv,\buv')=i\langle v,\nu'\rangle-i\langle\nu,v'\rangle\,,
\qquad \buv,\buv'\in\HH^2\,.
\equation(Gvvp)
$$

We are interested in the parameter-dependence of eigenvalues
that lie on the unit circle.
Let $s\mapsto A_s$ be a continuous curve of symplectic operators on $\HH^2$.
Let $\buv_s$ be eigenvector of $A_s$
with eigenvalue $\lambda_s$, both depending continuously on $s$.
If $\lambda_s$ lies on the unit circle and is simple, for some value $s=s_0$,
then the same is true for $s$ near $s_0$.
The reason is that, by symplecticity,
the spectrum of $A_s$ is invariant under
complex conjugation $z\mapsto\bar z$
and under inversion $z\mapsto z^{-1}$.
More can be said by using the the Krein signature of $\buv_s$,
which is defined to be the sign of
$$
G(\buv_s,\buv_s)=-2\,\Im \langle v_s,\nu_s\rangle\,.
\equation(Gvv)
$$
It is straightforward to check that $G(\buv_s,\buv_s)$ vanishes
unless $\lambda_s$ lies on the unit circle.
According to Krein theory, the only way that $\lambda_s$
can move off the unit circle, as $s$ is varied,
is for $\lambda_s$ to collide with an eigenvalue (for an eigenvector)
of opposite Krein signature.
This motivates the following

\claim Definition(SepSet)
Let $\Lambda=\circle\cup(0,\infty)$,
where $\circle$ denoted the unit circle in $\complex$.
Consider a finite subset $Z$ of $\circle\setminus\{1\}$
that contains at least two points.
This set defines a partition of $\Lambda\setminus Z$ into connected sets.
The set containing $1$ will be referred to as the ``cross''.
The other sets in this partition are subsets of $\circle$
and will be referred to as ``arcs''.
Given a symplectic operator $A$,
We say that $Z$ is a separating set for $A$,
if all eigenvalues of $A$
that lie in the same arc have the same (nonzero) Krein signature.
Furthermore, we impose that the cross contains
exactly two eigenvalues of $A$,
and that $Z$ contains no eigenvalues of $A$.

We note that the separating sets defined in [\rAKii]
were allowed to contain the point $1$.
But we only considered partitions of $\circle\setminus Z$,
since we did not allow bifurcations at $1$.
Here, the cross associated with $Z$ is needed
to control a pair of eigenvalues near $1$,
independently of whether these eigenvalues lie on $\circle$ or not.
Recall that the true system has an eigenvalue $1$,
and by symplecticity, this eigenvalue has an even multiplicity.
When considering finite-dimensional approximations,
this eigenvalue can split into multiple eigenvalues near $1$.

Using the above-mentioned fact about the Krein signature,
say in the form of Proposition 2.9 in [\rAKii],
we immediately obtain the following.

\claim Proposition(SepSetFam)
Let $Z$ be a finite subset of $\circle\setminus\{1\}$
that does not contain any eigenvalues of $A_s$ for any $s$.
Assume that one of the operators $A_s$
has the following property:
$Z$ is a separating set for $A_s$,
and all eigenvalues of $A_s$ belong to $\Lambda$.
Then each of the operators $A_s$ has this property.

Let $(s,t)\mapsto H^s(t)$ be a continuous family of
linear operators on $\HH$, indexed by $[-1,1]\times\real$.
Consider the flow on $\HH^2$ defined by the equation
$$
\partial_t\buv(t)=X_s(t)\buv(t)\,,\qquad
X_s(t)=\stwomat{0}{\id}{-H^s(t)}{0}\,.
\equation(dtvXv)
$$
Assume in addition that each operator $H^s(t)$ is self-adjoint.
Then a straightforward computation shows that
the time-$t$ maps $\Flow_s(t)$ for this flow are symplectic.
In what follows, we also assume that $H^s(t)$ is $2\pi$-periodic in $t$.
Then the map $\Flow_s(2\pi)$ is of particular interest.

The monotonicity property that we
mentioned earlier can be stated roughly as follows.
Assume that ${d\over ds}H^s(t)$ is positive for all $s$ and all $t$.
Then the eigenvalues of $\Flow_s(2\pi)$
that have negative (positive) Krein index
move (counter)clockwise on $\circle$, as $s$ is increased.

Formally, this follows from an explicit computation [\rAKii].
To make this statement more precise,
we need to avoid collisions of eigenvalues of opposite Krein signatures.
And for simplicity, we restrict now to affine families
$$
H^s=H^0+sD\,,\qquad s\in[-1,1]\,,
\equation(LinFamily)
$$
that are strongly increasing, in the sense that $D\gtstrong 0$.
Here, and in what follows,
if $C=C(t)$ and $D=D(t)$ are curves of self-adjoint linear operators on $\HH$,
then we define $D\gtstrong C$ or $C\ltstrong D$ to mean that
there exists $\eps>0$ such that
$D(t)-C(t)-\eps\id$ is a positive operator for all $t$.

An eigenvalue of $\Flow_s(2\pi)$ that lies on the unit circle
can be written as $\lambda=e^{2\pi i\eta}$.
The real number $\eta$ will be referred to as a Floquet number for $\Flow_s(2\pi)$.

\claim Proposition(Monotonicity)
{\rm(monotonicity)}
Assume that the family of operators $A_s=\Flow_s(2\pi)$
satisfies the hypotheses of \clm(SepSetFam).
Consider the eigenvalues of $A_s$ that lie in the arcs determined by $Z$.
Then the corresponding Floquet numbers $\eta_1,\eta_2,\ldots$
can be labeled in such a way that each $\eta_k$
is a real analytic function of the parameter $s$.
Furthermore, if the Krein signature of $\lambda_k$ is positive (negative)
then ${d\over ds}\eta_k$ is negative (positive).

This proposition is a consequence of \clm(SepSetFam),
and of Lemma 3.6 in [\rAKii].

\smallskip
Consider now a situation where $H^{-1}\ltstrong H\ltstrong H^1$,
as mentioned at the beginning of this subsection.
Since we can interpolate first between $H^{-1}$ and $H$,
and then between $H$ and $H^1$,
\clm(Monotonicity) suggests that each Floquet number for $\Flow(2\pi)$
can be bounded from above and below by the corresponding
Floquet numbers of $H^{-1}$ and $H^1$.
This is indeed the case, but the following suffices for our purpose.

\claim Proposition(FloquetSandwich)
Assume that $H^{-1}\ltstrong H\ltstrong H^1$.
Let $Z$ be a finite subset of $\circle\setminus\{1\}$
that does not contain any eigenvalues of $\Flow_s(2\pi)$, for any $s$.
Assume that, for some value of $s$, the operator $A=\Flow_s(2\pi)$
has the following property: $Z$ is a separating set for $A$,
and all eigenvalues of $A$ lie on $\Lambda$.
Then $A=\Flow(2\pi)$ has the same property.

The proof of this proposition is similar
to the proof of Corollary 3.8 in [\rAKii], so we omit it here.

\demo Remark(finitedim)
By \clm(ApproxEigen), it suffices to consider
the operators $H_k$ defined by \equ(Hn), if $k$ is chosen sufficiently large.
The results of this subsection will be applied with $H=H_k$.
Choosing $\HH=\HH_{0,n}$ with $n>k$,
the time-$2\pi$ map $\Flow(2\pi)$ leaves $\HH^2$ invariant,
and its spectrum does not depend on $n$.

\subsection Verifying separation

Motivated by \clm(SepSetFam), consider the task of verifying that $Z$
does not contain any eigenvalues of $A_s=\Flow_s(2\pi)$ for any $s$.
It is worth noting that this task simplifies
if we first verify that $Z$ is a separating set for $A_{-1}$,
and that all eigenvalues of $A_{-1}$ lie on $\Lambda$.
To see why, notice that every arc $\Gamma$ defined by $Z$
can be assigned a signature:
the signature of the eigenvalues of $A_{-1}$ that lie in $\Gamma$.
We may assume that $Z$ is ``minimal'',
in the sense that adjacent arcs have opposite signatures.
Consider now a point $z\in Z$,
and let $\Gamma$ be an arc that has $z$ in its boundary.
As $s$ is increased, starting from $-1$, the eigenvalues of $A_s$ in $\Gamma$
all move either toward $z$, or they all move away from $z$.
In the first case, we call $z$ a ``primary'' point of $Z$.
In the second case,
any eigenvalue that could possibly enter $\Gamma$ through $z$
must have the same signature as $\Gamma$,
so $z$ lies on the boundary of the cross.
Thus, if we verify that the eigenvalues of $A_s$
avoid all primary points of $Z$, as $s$ is increased from $-1$ to $1$,
and that $A_1$ has the same number of eigenvalues in the cross
as $A_{-1}$, then all points in $Z$ are being avoided.

\smallskip
We describe now a method for proving that
a given point $e^{2\pi i\eta}$ on the unit circle is not an eigenvalue
of any of the operators $\Flow_s(2\pi)$.
Consider first a fixed value of the parameter $s$.
Assume that $H^s$ is a self-adjoint linear operator on $\HH$
that depends continuously and $2\pi$-periodically on time $t$.
Let $\buv$ be an eigenvector of $\Flow_s(2\pi)$
with eigenvalue $\lambda=e^{2\pi i\eta}$,
and let $v$ be the first component of $\buv$.
Then the function $w=e^{-i\eta t}v$ is $2\pi$-periodic and satisfies the equation
$(\partial_t+i\eta)^2w=-H^s w$, or equivalently,
$$
M_s(\eta)w=0\,,\qquad M_s(\eta)=(\Kt+\eta)^2-H^s\,,\qquad
\Kt=-i\partial_t\,.
\equation(MsetaDef)
$$
Here $H^s w$ is defined pointwise by the equation $(H^s w)(t)=H^s(t)w(t)$.
If $\eta$ is real, then $M_s(\eta)$ is self-adjoint
as a linear operator on the Hilbert space $\buH=\rmL^2([0,2\pi],\HH)$
with the inner product
$$
\Langle w,w'\Rangle
={1\over\pi}\int_0^{2\pi}\bigl\langle w(t),w'(t)\bigr\rangle\,dt\,.
\equation(iiprod)
$$
To be more precise, the domain of $M_s(\eta)$
is the set of all functions $w\in\buH$ with the property
that $\Kt w$ belongs to $\buH$.
Then $M_s(\eta)$ is a Fredholm operator on $\buH$,
and in particular, the spectrum of $M_s(\eta)$ consists
of isolated eigenvalues with finite multiplicity.
It is straightforward to show that the eigenfunctions
of $M_s(\eta)$ are continuous,
and that any nonzero vector $w$ in the null space of $M_s(\eta)$
yields an eigenvalue $\buv$ of $\Flow_s(2\pi)$
with eigenvalue $\lambda=e^{2\pi i\eta}$.
For details we refer to [\rAKii].

We need the above only for operators $H^s$ in an affine family
\equ(LinFamily), where $(Dw)_j=d_jw_j$ for
a bounded sequence $j\mapsto d_j$ of positive real numbers.
Consider the affine family $H^s=H^0+sD$ for $s\in[-1,1]$.
With $M_s(\eta)$ as defined in \equ(MsetaDef),
our goal is to show that, for some given $\eta\in\real$,
none of the operators $M_s(\eta)$ with $s\in[-1,1]$ has an eigenvalue zero.

It is convenient to replace $M_s(\eta)$
by a bounded linear operator $\hat M_s(\eta)$ as follows.

\claim Definition(theta)
For every integer $k$, define $\theta_k=\max(1,|k|)^{-1}$.
Denote by $\theta$ the (unique) continuous linear operator
on $\buH$ with the property that for each $j$,
if $w_j(t)=e^{ikt}$ for all $t$, then $(\theta w)_j(t)=\theta_ke^{ikt}$ for all $t$.
If $M$ is any linear operator on $\buH$, then we define
$\hat Mw=\theta M\theta w$, whenever $w\in\buH$
and $\theta w$ belongs to the domain of $M$.

Clearly, $M_s(\eta)$ has an eigenvalue zero if and only if
$\hat M_s(\eta)$ has an eigenvalue zero.
The operators $\theta(\Kt+\eta)^2\theta$ and $\hat D$ are trivial to represent.
The operator $\hat H^0$ is less easy to handle.
But it is compact, so we approximate it by a simpler (finite rank)
operator $\check H^0$.

In order to show that $e^{2\pi i\eta}$ is not an eigenvalue of
$\Flow_s(2\pi)$, it suffices now to verify the hypotheses of the following lemma.

\claim Proposition(SepCond) {\rm[\rAKii]}
Consider parameters values $-1=s_0<s_1<\ldots<s_m=1$.
Let $C>\bigl\|\hat H^0-\check H^0\bigr\|$.
Assume that the operator
$\check M_{s_j}(\eta)=\theta(\Kt+\eta)^2\theta-\check H^0-s\hat D$
has no eigenvalues in $[-C,C]$,
and that
$$
(s_j-s_{j-1})\bigl\|\hat D\bigr\|<2C\,,
\equation(SepCond)
$$
for $j=1,2,\ldots,m$.
Then none of the operators $\hat M_s(\eta)$
has an eigenvalue zero.

\subsection Proof of \clm(stability)

Before we can apply \clm(FloquetSandwich),
we need to find useful upper and lower bounds
of the form $H^0-D\ltstrong H_n\ltstrong H^0+D$
on the operators $H_n$ defined by the equation \equ(Hn).
To be more precise, such bounds are needed only for sufficiently large $n$.
Thus, let first determine upper and lower bounds
on the full operator $H$, given by \equ(HqDef).
Instead of $H(q(t))$,  we write here $H(t)$ or just $H$.
We assume that $\beta_j\ge 0$ for all $j$.

The space considered here is $\HH=\ell^2(\integer)$.
We start by defining a self-adjoint truncation $H^o$ of the operator $H$.
After fixing a cutoff $j_\ast>0$, $H^o$ is defined by the quadratic form
$$
\langle v,H^o v\rangle
=\sum_{\sigma-j_\ast\le j<j_\ast}
\alpha_{j-\sigma+1}|v_{j+1}-v_j|^2
+\sum_{\sigma-j_\ast\le j\le j_\ast}\bigl(\beta_j-\bar\beta\bigr)|v_j|^2
+\bar\beta\langle v,v\rangle\,,
\equation(vHov)
$$
where $\alpha$ and $\beta$ are the functions defined in \equ(abDef),
and where $\bar\beta=\omega^{-2}\psi_2$.
The truncation error $\EE=H-H^o$ is then given by
$$
\langle v,\EE v\rangle
=\!\!\sum_{j\ge j_\ast{\rm~or~}j<\sigma-j_\ast}\!\!\alpha_{j-\sigma+1}|v_{j+1}-v_j|^2
+\!\!\sum_{j>j_\ast{\rm~or~}j<\sigma-j_\ast}\!\!\bigl(\beta_j-\bar\beta\bigr)|v_j|^2\,.
\equation(vEEv)
$$
In order to estimate $\EE$,
we determine for $j=j_\ast$ and for $j=\sigma-j_\ast$
an interval $[\alpha_j^{-},\alpha_j^{+}]$ that includes $\{0,\alpha_j\}$.
In addition, we determine an interval $[\alpha_\infty^{-},\alpha_\infty^{+}]$
that includes $\{0,\alpha_j\}$ whenever $j<1-j_\ast$ or $j>j_\ast-\sigma$.
And we choose constants $\gamma^{-}\le 0\le\gamma^{+}$ such that
$$
\gamma^{-}\le 4\alpha_\infty^{\pm}+\beta_j-\bar\beta\le\gamma_j^{+}\,,\qquad
\hbox{if $j<\sigma-j_\ast$ or $j>j_\ast$}\,.
\equation(tailbounds)
$$
Now define
$$
\pm d_j^{\pm}=\cases{2\alpha_{j-s+1}^{\pm} &if $j=j_\ast$ ,\cr
            2\alpha_{j-s}^{\pm} &if $j=\sigma-j_\ast$ ,\cr
            \gamma^{\pm} &if $j<\sigma-j_\ast$ or $j>j_\ast$ ,\cr
            0 &if $\sigma-j_\ast<j<j_\ast$ .\cr}
\equation(djDef)
$$
Using the trivial inequality $|x-y|^2\le 2|x|^2+2|y|^2$, we find that
$$
-D^{-}\ltstrong\EE\ltstrong D^{+}\,,\qquad
D^{\pm}=\diag\bigl(d^{\pm}+\epsilon)\,,
\equation(EESandwich)
$$
for any $\epsilon>0$. Finally, define
$$
H^s=H^0+sD\,,\qquad
H^0=H^o+\thalf\bigl(D^{+}-D^{-}\bigr)\,,\qquad
D=\thalf\bigl(D^{+}+D^{-}\bigr)\,.
\equation(HoD)
$$
Then $D\gtstrong 0$ and
$$
H^{-1}=H^o-D^{-}\ltstrong H\ltstrong H^o+D^{+}=H^1\,.
\equation(HoDmDpSandwich)
$$
It is straightforward to check that the same
holds if $H$ is replaced by any of the operators $H_n$
with $n$ sufficiently large.

\demo Remark(cleanHo)
The function $H^o$ that is used in our computer-assisted proof
differs from \equ(vHov) in the sense that $\alpha_j$ and $\beta_j$
are replaced by function $\alpha_j^o$ and $\beta_j^o$
that are very close to $\alpha_j$ and $\beta_j$, respectively.
The value of $\epsilon>0$ in the definition \equ(EESandwich)
is chosen to (over)compensate for the resulting error.

Notice that a chain $v$ that is supported at a single point
$j<\sigma-j_\ast$ or $j>j_\ast$ is an eigenvector of $H^s$,
with eigenvalue
$$
\textstyle
\mu_s=\bar\beta+{1-s\over 2}(\gamma^{-}-\epsilon)
+{1+s\over 2}(\gamma^{+}+\epsilon)\,.
\equation(highEigen)
$$
The corresponding eigenvalues of $\Flow_s(2\pi)$ need to be considered as well
in our application of \clm(FloquetSandwich).
In $\HH$ they have infinite multiplicity,
but when considering $H_n$ in place of $H$,
only the eigenvalues with eigenvectors in $P_n\HH$ are relevant.
In order to compute their Krein signature, write $\mu_s=-\eta_s^2$.
Then the corresponding eigenvalues for $X_s$ are $\pm i\eta_s$.
So we need $\mu_s\le 0$ in order for $\Flow_s(2\pi)$ to be spectrally stable.
Assume that $\eta_s>0$.
Using \equ(Gvv), one easily finds that the eigenvector
for the eigenvalue $e^{2\pi i\eta_s}$ of $\Flow_s(2\pi)$ has a negative Krein signature,
while $e^{-2\pi i\eta_s}$ has a positive Krein signature.

The following two lemmas are proved with the assistance
of a computer, as described in Section 5.

\claim Lemma(separation)
For each of the solutions $2$, $3$, $4$, $8$, and $13$,
there exists a family of operators $s\mapsto H^s$ as described above,
satisfying $H^{-1}\ltstrong H_n\ltstrong H^1$ for large $n$.
In addition, there exists a common separating set $Z$ for both $\Flow_{\pm 1}(2\pi)$,
a finite rank operator $\check H^0$, a constant $C>0$,
and parameter values $-1=s_0<s_1<\ldots<s_m=1$,
such that the hypotheses of \clm(SepCond) are satisfied,
for every primary point $e^{2\pi i\eta}$ in $Z$.

For the definition of a primary point in $Z$,
we refer to Subsection 4.3.

\claim Lemma(unstable)
For each of the solutions $1$ and $11$,
there exist real numbers $\alpha_\ast^{\pm}$, $\beta_\ast^{\pm}$,
and $\mu>1$ such that the hypotheses of \clm(unstableBound)
are satisfied, with $H^o$ as described above.

\medskip
We note that the separating set $Z$ described in \clm(separation)
is determined by computing accurate bounds
on the eigenvalues of the time-$2\pi$ maps $\Flow_{\pm 1}(2\pi)$
associated with the operator $H^{\pm 1}$.
The nontrivial part of $\Flow_{\pm 1}(2\pi)$ is just a $2k_\ast\times 2k_\ast$
matrix, where $k_\ast=2j_\ast+1-\sigma$.
It is obtained by integrating the flow $\dot\buv=X_{\pm 1}\buv$
associated with the second order equation $\ddot v=-H^{\pm 1}v$.

In order to make this part of our programs [\rFiles] more transparent,
let us write down the equations that are being integrated.
To simplify notation, consider the operator $H^o$ in place of $H^{\pm 1}$.
After a change of variables (indices) $g_k=v_{k-j_\ast-1+\sigma}$,
the equation $\ddot v=-H^o v$ becomes
$$
\ddot g_k
=\alpha_{k-j_\ast}g_{k+1}
-(\alpha_{k-j_\ast-1}+\alpha_{k-j_\ast}+\beta_{k-j_\ast-1+s})g_k
+\alpha_{k-j_\ast-1}g_{k-1}\,,
\equation(ddotgCenter)
$$
for $1<k<k_\ast$, and
$$
\eqalign{
\ddot g_1
&=\alpha_{1-j_\ast}g_2
-(\alpha_{1-j_\ast}+\beta_{-j_\ast+s})g_1\,,\cr
\ddot g_{k_\ast}
&=-(\alpha_{j_\ast-s}+\beta_{j_\ast})g_{k_\ast}
+\alpha_{j_\ast-s}g_{k_\ast-1}\,.\cr}
\equation(ddotGTail)
$$

\medskip
Based on Lemmas \clmno(separation) and \clmno(unstable),
we can now give a

\proofof(stability).
Consider first one of solutions $2$, $3$, $4$, $8$, and $13$,
that we claim to be spectrally stable.
By \clm(separation), we can apply \clm(SepCond) to conclude
that none of the points in $Z$
is an eigenvalue of any of the operators $\Flow_s(2\pi)$.
Here, we have also used the argument
given at the beginning of Subsection 4.3,
which shows that it suffices to check the primary points.
Now we can apply \clm(FloquetSandwich),
with $H_n$ in place of $H$, for $n$ sufficiently large.
It shows that $Z$ is a separating set for $\Flow_n(2\pi)$,
and that all eigenvalues of $\Flow_n(2\pi)$ lie on $\Lambda$.
Taking $n\to\infty$ along a suitable subsequence,
we conclude from \clm(LimEigen) that all eigenvalues
of $\Flow(2\pi)$ lie on the unit circle
and are bounded away from $1$,
with the possible exception of two eigenvalues
on the closure of the cross determined by $Z$.
But we already know that $\Flow(2\pi)$ has an eigenvalue $1$,
as mentioned in \dem(translation),
and this eigenvalue must have an even multiplicity by symplecticity.
This implies that all eigenvalues of $\Flow(2\pi)$
lie on the unit circle.

Next, consider one of solutions
that we claim to be spectrally unstable.
As mentioned at the beginning of Subsection 4.1,
it suffices to consider the solutions $1$ and $11$.
In these cases, \clm(unstable) and \clm(unstableBound)
imply that $\Flow(2\pi)$ has at least one real eigenvalue larger than $1$.
This concludes the proof of \clm(stability).
\qed

\section Computer estimates

In order to complete our proof of
Theorems \clmno(breathers), \clmno(stability), and \clmno(multi),
we need to verify the assumptions of the
Lemmas \clmno(contraction), \clmno(bumpcontr), \clmno(separation), and \clmno(unstable).
The strategy is to reduce each of these lemmas
to successively simpler propositions,
until the claims are trivial numerical statements that can be
(and have been) verified by a computer.
This part of the proof is written in the programming language Ada [\rAda]
and can be found in [\rFiles].

The following is meant to be a rough guide for the reader
who wishes to check the correctness of our programs.
The first part of the above-mentioned reduction
is organized by the main program {\tt Run\_All}.
It divides the given task among five standalone procedures.
The first is {\tt Approx\_Fixpt}, which is purely numerical
and computes the finite-rank part $A'=A-\id$ of the operator $A$
that appears in \equ(contr).
The approximate solution $\bar q$ is read from the {\tt data} directory,
and the necessary parameters are specified in the Ada package {\tt Params}.
(If desired, {\tt Approx\_Fixpt} can be used also
to improve the quality of the approximate solution.)
Now that the map $\FF$ is well-defined,
the procedure {\tt Check\_Fixpt} is called to
verify the assumptions of \clm(contraction) and \clm(bumpcontr).
At this point, we have an enclosure for the fixed point $q$ of $G$.
Enclosures for chains in $\BBB$ are represented by the data type {\tt FChain},
using enclosures of type {\tt CosSin1} for functions in $\AA_r^\pars$.
Data associated with $q$ that are needed later,
such as the functions $\alpha$, $\beta$,
and upper bounds on the numbers $d_j^{\pm}$ defined in \equ(djDef),
are computed and saved by the procedure {\tt Save\_Data}.
This procedure also determines a bound
{\tt NPD} on the operator norm $\|H-H^o\|$ that appears in \equ(uBoundOne).
Bounds on the maps $\Flow_{\pm 1}(2\pi)$ and on its eigenvalues are determined
by the procedures {\tt Phi2Pi} and {\tt Eigen}.
For the solutions that are expected to be unstable,
{\tt Eigen} also calls the procedure {\tt ScalVectors.Phi.Check\_Unstable}
to verify the the assumptions of \clm(unstable).
For the solutions that are expected to be stable,
{\tt Run\_All} calls the procedure {\tt Separation}
to verify the assumptions of \clm(separation).

The next steps in the reduction process require specialized knowledge
and tools, so each of the above-mentioned procedures first
instantiates a few specialized Ada packages
and then hands the task to some procedure(s) that are implemented in those packages.
An Ada package is simply a collection of definitions and procedures,
centered around a few specific data types.
In particular, the package {\tt CosSins1} and its child {\tt CosSins1.Chain}
implement basic bounds involving the data types
{\tt CosSin1} and {\tt FChain}, respectively.
The type {\tt CosSin1} is equivalent to the type {\tt Fourier}
that is used and documented in [\rAGT].
Our type {\tt FChain} is in essence an array of {\tt CosSin1},
indexed by $\JJ=\{j\in\integer:\sigma-\ell\le j\le\ell\}$.
If {\tt Q} is an {\tt FChain} specifying an enclosure for a chain $q\in\BBB$,
then the components ${\tt Q}_{\sigma-\ell}$ and ${\tt Q}_\ell$
consist of error bounds on the tails $(\ldots,q_{\sigma-\ell-1},q_{\sigma-\ell})$
and $(q_\ell,q_{\ell+1},\ldots)$, respectively.
The remaining components ${\tt Q}_j$ define enclosures for the functions $q_j\in\AA_r^\pars$,
with $\sigma-\ell<j<\ell$. In our programs, $\ell$ is named {\tt JEMax}.
And the cutoff $j_\ast<\ell$ used in \equ(vHov) is named {\tt JAst}.

As can be seen in {\tt Check\_Fixpt},
the specialized bounds that are needed in the proofs of \clm(contraction) and \clm(bumpcontr)
are implemented in the child package {\tt CosSins1.Chain.Fix}.
This includes bounds {\tt GMap} and {\tt DGMap} on the maps $G$ and $DG$, respectively.
Similarly, the proof of \clm(separation)
is organized by the procedure {\tt CheckEta}
in the package {\tt CosSins1.Chain.Pairs.FlokM}.
As the package structure indicates,
bounds defined in {\tt CosSins1.Chain.Pairs.FlokM}
are reduced in stages to bounds defined in {\tt CosSins1},
and those reduce further to bounds on data of type {\tt Scalar}, etc.
Following these instructions,
a computer ends up with a finite number of basic numerical operations,
which are carried out with rigorous upper and (if necessary) lower bounds.

All this is described in full detail by the source code of our programs [\rFiles].
But some remarks may be in order concerning the choice of algorithms.
Whenever an implicit equation needs to be solved,
our approach is the same as for the equation $G(q)=q$.
After determining an approximate solution $\bar q$,
we use the contraction mapping theorem for a Newton-type map $\FF$
to obtain a rigorous bound on the error $q-\bar q$.
This approach is used e.g.~to obtain bounds on the eigenvalues $\lambda_k$
of the symplectic matrix for the nontrivial part of $\Flow_{\pm 1}(2\pi)$ or $\Flow_o(2\pi)$,
after determining a polynomial
whose roots are the numbers $\half\lambda_k+\half\lambda_k^{-1}$.
The computation of the matrix itself is entirely explicit:
here we use a Taylor method to integrate the nontrivial part
of the vector field $X_{\pm 1}$ (for \clm(separation)) or $X_o$ (for \clm(unstable))
associated with the operators $H^o\pm D_{\pm 1}$ or $H^o$, respectively,
described in Subsection 4.4. See also the comments after \clm(unstable).
Verifying the assumptions of \clm(separation) is an explicit computation as well.
Here, ``computing'' an object means finding a rigorous enclosure
(specified by finitely many representable numbers) for that object.
To prove that the operator $L=\check M_{s_j}(\eta)$ described in \clm(SepCond)
has no eigenvalue in $[-C,C]$, we simply compute the inverse of $L$ and
check that $\|L^{-n}\|<C^n$ for some positive integer $n$ (a power of $2$).
A more detailed description of the algorithms
used to integrate a vector field and to compute eigenvalues can be found in [\rAKii],
where we considered a similar spectral problem.

We will not explain here the more basic ideas and techniques
underlying computer-assisted proofs in analysis.
This has been done to various degrees in many other papers, including [\rAKT,\rAKii].
As far as our proof of the
Lemmas \clmno(contraction), \clmno(bumpcontr), \clmno(separation), and \clmno(unstable)
is concerned, the ultimate reference is the source code of our programs [\rFiles].
For the set of representable numbers ({\tt Rep})
we choose either standard [\rIEEE] extended floating-point numbers (type {\tt LLFloat})
or high precision [\rMPFR] floating-point numbers (type {\tt MPFloat}),
depending on the precision needed. Both types support controlled rounding.
Our programs were run successfully on a standard
desktop machine, using a public version of the gcc/gnat compiler [\rGnat].
Instructions on how to compile and run these programs
can be found in the file README that is included with the source code [\rFiles].

\references

{\ninepoint

\item{[\rKato]} T.~Kato,
{\sl Perturbation Theory for Linear Operators},
Springer Verlag, 1976.

\item{[\rAGb]} G.~Arioli, F.~Gazzola,
{\sl Existence and numerical approximation of periodic motions of an infinite lattice of particles},
ZAMP {\bf 46}, 898--912 (1995).

\item{[\rAGa]} G.~Arioli, F.~Gazzola,
{\sl Periodic motions of an infinite lattice of particles with nearest neighbor interaction},
Nonlin. Anal. TMA {\bf 26}, 1103--1114 (1996).

\item{[\rAGT]} G.~Arioli, F.~Gazzola, S.~Terracini,
{\sl Multibump periodic motions of an infinite lattice of particles},
Math. Zeit. {\bf 223}, 627--642 (1996).

\item{[\rRabi]} P.H.~Rabinowitz,
{\sl Multibump solutions of differential equations: an overview},
Chinese J. Math. {\bf 24}, 1--36 (1996).

\item{[\rBK]} O.M.~Braun, Y.S.~Kivshar,
{\sl The Frenkel-Kontorova model},
Springer Verlag, Berlin (2004).

\item{[\rAKT]} G.~Arioli, H.~Koch, S.~Terracini,
{\sl Two novel methods and multi-mode periodic solutions
for the Fermi-Pasta-Ulam model},
Commun. Math. Phys. {\bf 255}, 1--19 (2004).

\item{[\rBerIzr]} G.P.~Berman, F.M.~Izraileva,
{\sl The Fermi-Pasta-Ulam problem: Fifty years of progress},
Chaos {\bf 15}, 015104 (2005).

\item{[\rAPankov]} A.~Pankov,
{\sl Traveling Waves And Periodic Oscillations in Fermi-Pasta-Ulam Lattices},
Imperial College Press (2005).

\item{[\rGalla]} G.~Gallavotti (editor),
{\sl The Fermi-Pasta-Ulam problem. A status report},
Lecture Notes in Physics {\bf 728}. Springer, Berlin, Heidelberg (2008).

\item{[\rGFii]} A.V.~Gorbach, S.~Flach,
{\sl Discrete breathers --- Advances in theory and applications},
Phys. Reports {\bf 467}, 1--116 (2008).

\item{[\rFLMii]} E.~Fontich, R.~de la Llave, P.~Martin,
{\sl Dynamical systems on lattices with decaying interaction II:
Hyperbolic sets and their invariant manifolds},
J. Differ. Equations {\bf 250}, 2887--2926 (2011).

\item{[\rPS]} D.~Pelinovsky, A.~Sakovich,
{\sl Multi-site breathers in Klein-Gordon lattices:
stability, resonances, and bifurcations},
Nonlinearity {\bf 25}, 3423--3451 (2012).

\item{[\rAKii]} G.~Arioli, H.~Koch,
{\sl Spectral stability for the wave equation with periodic forcing},
J. Differ. Equations {\bf 265}, 2470--2501 (2018).

\item{[\rFiles]} G.~Arioli, H.~Koch.
The source code for our programs, and data files, are available at\hfill\break
\pdfclink{0 0 1}{{\tt http://web.ma.utexas.edu/users/koch/papers/breathers/}}
{http://web.ma.utexas.edu/users/koch/papers/breathers/}

\item{[\rAda]} Ada Reference Manual, ISO/IEC 8652:2012(E),
available e.g. at\hfil\break
\pdfclink{0 0 1}{{\tt http://www.ada-auth.org/arm.html}}
{http://www.ada-auth.org/arm.html}

\item{[\rGnat]}
A free-software compiler for the Ada programming language,
which is part of the GNU Compiler Collection; see
\pdfclink{0 0 1}{{\tt http://gnu.org/software/gnat/}}{http://gnu.org/software/gnat/}

\item{[\rIEEE]} The Institute of Electrical and Electronics Engineers, Inc.,
{\sl IEEE Standard for Binary Float\-ing--Point Arithmetic},
ANSI/IEEE Std 754--2008.

\item{[\rMPFR]} The MPFR library for multiple-precision floating-point computations
with correct rounding; see
\pdfclink{0 0 1}{{\tt http://www.mpfr.org/}}{http://www.mpfr.org/}

}

\bye